\documentclass[aap,MSNbibl,citesort,dvips]{arximspdf}
\usepackage[mathcal]{euscript}
\usepackage{graphicx}

\doi{10.1214/11-AAP800}
\volume{22}
\issue{4}
\pubyear{2012}
\firstpage{1495}
\lastpage{1540}

\makeatletter

\newcommand{\sfPi}{\mbox{{$\mathsf{I}\!\rule[6.2pt]{5pt}{0.9pt}\!\mathsf{I}$}}}

\newcommand{\cdv}{\cdot\,|}

\newcommand{\nompageref}{\pageref}

\newcommand{\xrightarrow}[1]{\stackrel{#1}{\hbox to 1cm{\rightarrowfill}}}

\newtheorem{theorem}{Theorem}[section]
\newtheorem{cor}[theorem]{Corollary}
\newtheorem{lem}[theorem]{Lemma}
\newtheorem{prop}[theorem]{Proposition}

\newproclaim{aspt}[theorem]{Assumption}
\newproclaim{rem}[theorem]{Remark}
\newproclaim{example}[theorem]{Example}

\makeatother

\begin{document}
\begin{frontmatter}

\title{Ergodicity and stability of the conditional distributions of
nondegenerate Markov chains\thanksref{T1}}
\runtitle{Ergodicity and stability of conditional distributions}

\thankstext{T1}{Supported in part by NSF Grant DMS-10-05575.}

\begin{aug}
\author[A]{\fnms{Xin Thomson} \snm{Tong}\ead[label=e2]{xintong@princeton.edu}}
\and
\author[A]{\fnms{Ramon} \snm{van Handel}\corref{}\ead[label=e3]{rvan@princeton.edu}}
\runauthor{X. T. Tong and R. van Handel}
\affiliation{Princeton University}
\address[A]{Sherrerd Hall\\
Princeton University \\
Princeton, New Jersey 08544 \\
USA \\
\printead{e2} \\
\hphantom{E-mail: }\printead*{e3}} %adresu isvedimo komanda gale!
\end{aug}

% HISTORY:
\received{\smonth{1} \syear{2011}}
\revised{\smonth{8} \syear{2011}}

% ABSTRACT
%
\begin{abstract}
We consider a bivariate stationary Markov chain $(X_n,Y_n)_{n\ge0}$ in
a Polish state space, where only the process $(Y_n)_{n\ge0}$ is
presumed to be observable. The goal of this paper is to investigate the
ergodic theory and stability properties of the measure-valued process
$(\Pi_n)_{n\ge0}$, where $\Pi_n$ is the conditional distribution of
$X_n$ given $Y_0,\ldots,Y_n$. We show that the ergodic and stability
properties of $(\Pi_n)_{n\ge0}$ are inherited from the ergodicity of
the unobserved process $(X_n)_{n\ge0}$ provided that the Markov chain
$(X_n,Y_n)_{n\ge0}$ is nondegenerate, that is, its transition kernel
is equivalent to the product of independent transition kernels. Our
main results generalize, subsume and in some cases correct previous
results on the ergodic theory of nonlinear filters.
\end{abstract}

% KEYWORDS
%

\end{frontmatter}

%s1 #&#
\section{Introduction}\label{intro}

In this paper we will consider a bivariate Markov chain $(X_n,Y_n)_{n\ge
0}$ taking values in a Polish state space. Only the process $(Y_n)_{n\ge
0}$ is presumed to be directly observable to us, and we aim to estimate
the state~$X_n$ of the unobserved process given the observed data
$Y_0,\ldots,Y_n$ to date. This is the quintessential setup in problems
with partial information, and models of this type can therefore be found
in a wide range of applications~\cite{CMR05}.

We will be concerned, in particular, with the ergodic theory and stability
properties of the measure-valued process $(\Pi_n)_{n\ge0}$ defined by the
conditional distributions
$\Pi_n=\mathbf{P}(X_n\in\cdv Y_0,\ldots,Y_n)$, which is called the
nonlinear filter. It is not difficult to show that, in general, the
processes $(\Pi_n,Y_n)_{n\ge0}$ as well as $(\Pi_n,X_n,Y_n)_{n\ge0}$ are
themselves Markovian, and a typical question that we will aim to answer is
whether ergodicity of the underlying model $(X_n,Y_n)_{n\ge0}$ implies
ergodicity of\vadjust{\goodbreak} the extended Markov chain $(\Pi_n,X_n,Y_n)_{n\ge0}$ in a
suitable sense. Questions of this type date back at least to the work of
Blackwell~\cite{Bla57} and Kunita~\cite{Kun71}. Beside the intrinsic
probabilistic interest in the development of a \textit{conditional}
ergodic theory of Markov chains, ergodicity of the filter has substantial
practical relevance to understanding the performance of nonlinear
filtering and its numerical approximations over a~long time horizon;
cf.
\cite{Kun71,BK01,vH09spa}, and see~\cite{vH10,CR11} for further
references.

%f1 #&#
\begin{figure}
\begin{tabular}{@{}c@{}}

\includegraphics{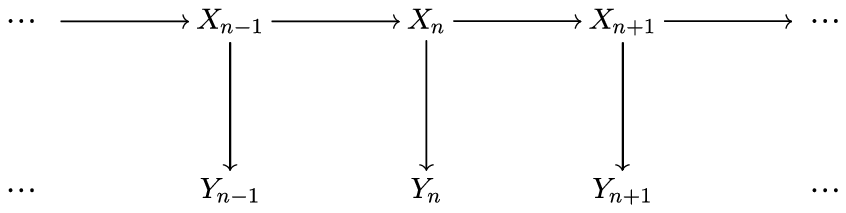}
\\
(a)\\[4pt]

\includegraphics{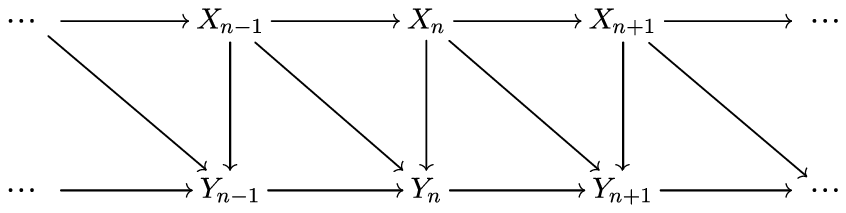}
\\
(b)\\[4pt]

\includegraphics{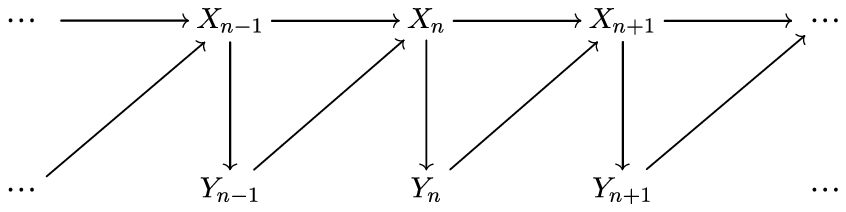}
\\
(c)\\[4pt]

\includegraphics{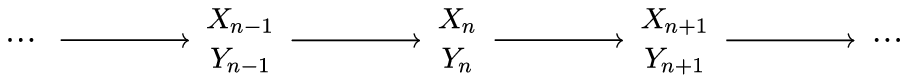}
\\
(d)
\end{tabular}
\caption{Dependence structure of \textup{(a)} a classical hidden Markov
model \protect\cite{CMR05}; \textup{(b)} a generalized hidden Markov
model \protect\cite{DS05,DMR04}; \textup{(c)} a hidden Markov model with
correlated noise \protect\cite{Bud02}; \textup{(d)}
a general Markov model.}\label{fighmm}
\end{figure}

Much of the literature on the topic of this paper is concerned with the
setting of a classical hidden Markov model whose dependence structure
is illustrated in Figure~\ref{fighmm}(a); here the unobserved process
$(X_n)_{n\ge0}$ is assumed to be itself Markovian and the observations
$(Y_n)_{n\ge0}$ are conditionally
independent.\setcounter{footnote}{1}\footnote{The continuous time version
of this model, known as a
Markov additive process, is also widely studied in the
literature in various special cases (such as white noise
or counting observations; see~\cite{Yor77} for a unified
view). We have restricted ourselves in this paper to
discrete time models for simplicity. All our results
are easily extended to the continuous time setting as in
\cite{vH10}, Section 6.}
In this special case $(\Pi_n)_{n\ge0}$ is also Markovian, and two basic
questions have been considered.
\begin{longlist}[(1)]
\item[(1)] Does $(\Pi_n)_{n\ge0}$ possess a unique invariant measure,
assuming $(X_n)_{n\ge0}$ does?
\end{longlist}
For the second question, let $\tilde{\mathbf P}$ and $\mathbf{P}$ be the
laws of the Markov chain $(X_n,Y_n)_{n\ge0}$ with initial laws
$\tilde{\mathbf P}(X_0\in\cdot)\ll\mathbf{P}(X_0\in\cdot)$, and
let $\tilde\Pi_n=\tilde{\mathbf P}(X_n\in\cdv Y_0,\ldots,Y_n)$.
\begin{longlist}[(2)]
\item[(2)] Is $(\Pi_n)_{n\ge0}$ asymptotically stable in the
sense that $|\tilde\Pi_n(f)-\break\Pi_n(f)|\xrightarrow{n\to\infty}0$
in $\tilde{\mathbf P}$-probability for every bounded continuous
function $f$?
\end{longlist}
These and related questions were studied in great generality by
Kunita~\cite{Kun71,Kun91}, Stettner~\cite{Ste89}, and Ocone and Pardoux
\cite{OP96} (see~\cite{Bud03,BCL04,vH10,CR11} for further references).
Kunita and Stettner state that the answer to the first question is affirmative
provided that the stationary process $(X_n)_{n\in\mathbb{Z}}$ is purely
nondeterministic, that is,
\[
\bigcap_{n\ge0}\mathcal{F}^X_{-n} \mbox{ is }
\mathbf{P}\mbox{-trivial},
\]
where $\mathbf{P}$ is the stationary law of the two-sided process
$(X_n,Y_n)_{n\in\mathbb{Z}}$ and $\mathcal{F}^X_n =
\sigma\{X_k\dvtx{-\infty}< k\le n\}$. Ocone and Pardoux state that the
answer to the second question is affirmative under the same assumption.
Unfortunately, the proofs of these results contain a serious error,
as was pointed out by Baxendale, Chigansky and Liptser~\cite{BCL04}.
Indeed, the crucial step in the proofs is the identity
\[
\bigcap_{n\ge0} \mathcal{F}^Y_0\vee\mathcal{F}^X_{-n}
\mathop{\stackrel{?}{=}}
\mathcal{F}^Y_0, \qquad\mathbf{P}\mbox{-a.s.},
\]
where $\mathcal{F}^Y_0=\sigma\{Y_k\dvtx{-\infty}<k\le0\}$. It is tempting
to exchange the order of the intersection $\cap$ and supremum $\vee$ of
$\sigma$-fields, which would allow us to conclude this identity from the
assumption that $(X_n)_{n\in\mathbb{Z}}$ is purely nondeterministic. But
such an exchange cannot be taken for granted (see~\cite{CY03}, page 30)
and requires proof. In the filtering setting, various counterexamples
given in~\cite{BCL04,vH11} show that the answers to the above questions
may indeed be negative even when $(X_n)_{n\in\mathbb{Z}}$ is purely
nondeterministic, in contradiction with the conclusions of
\cite{Kun71,Kun91,Ste89,OP96}.

Before we proceed, let us briefly recall a simple counterexample from
\cite{BCL04,vH11} that will be helpful in understanding the problems
addressed in this paper.
%
%ex1.1 #&#
\begin{example}
\label{exhmm}
Let $(\xi_n)_{n\in\mathbb{Z}}$ be an i.i.d. sequence of random
variables taking the values $\{0,1\}$ with equal probability
under $\mathbf{P}$, and define
\[
X_n = (\xi_{n},\xi_{n+1}),\qquad
Y_n = |\xi_{n+1}-\xi_{n}|.
\]
Then $(X_n)_{n\in\mathbb{Z}}$ is an ergodic Markov chain in
$\{00,01,10,11\}$ that is purely nondeterministic by the
Kolmogorov zero--one law, and $(X_n,Y_n)_{n\ge0}$ is
a~hidden Markov model as in Figure~\ref{fighmm}(a).\vadjust{\goodbreak}
But clearly $\xi_n = (\xi_1+Y_1+\cdots+Y_{n-1})\operatorname{mod}2$,
so that
\[
\tilde\Pi_n(f) = f(X_n),\qquad
\Pi_n(f) = \frac{f(X_n) + f(11-X_n)}{2},\qquad
\tilde{\mathbf P}\mbox{-a.s.}
\]
where we defined $\tilde{\mathbf P}( \cdot)=\mathbf{P}( \cdv X_0=00)$.
Thus the filter is not asymptotically stable, and one may similarly
establish that it admits distinct invariant measures.
\end{example}

One feature of the model of Example~\ref{exhmm} is that it possesses
degenerate observations in the sense that $Y_n$ is a function of $X_n$
without any additional noise. The phenomenon illustrated here turns out
to disappear when some independent noise is added to the observations, for
example, $Y_n = |\xi_{n+1}-\xi_n|+\eta_n$ where $(\eta_n)_{n\in\mathbb{Z}}$
is an i.i.d. sequence such that the law of $\eta_0$ has a~nowhere vanishing
density. In~\cite{vH10}, one of the authors developed this idea to
establish ergodicity and stability properties of the nonlinear
filter under very general assumptions. To this end, let
$(X_n,Y_n)_{n\in\mathbb{Z}}$ be a stationary hidden Markov model under
$\mathbf{P}$, and assume that:
\begin{longlist}[(2)]
\item[(1)] $(X_n)_{n\in\mathbb{Z}}$ is absolutely regular:
\mbox{$\mathbf{E}(\|\mathbf{P}(X_n\in\cdv X_0)-
\mathbf{P}(X_n\in\cdot)\|_{\mathrm{TV}})\to0$}.
\item[(2)] The observations are nondegenerate:
$\mathbf{P}(Y_n\in A|X_n) = \int_A g(X_n,y) \varphi(dy)$
for some strictly positive density $g(x,y)>0$ and reference measure
$\varphi$.
\end{longlist}
Then the above exchange of intersection and supremum of $\sigma$-fields
is permitted, and the filter is stable~\cite{vH10} and uniquely ergodic
\cite{vH11}. Intuitively, nondegeneracy (which formalizes the notion of
``noisy'' observations) rules out the singular observation structure
that causes the exchange of intersection and supremum to fail
in Example~\ref{exhmm}. However, this intuition should not be
taken too literally, as a more difficult example in~\cite{vH11} shows
that the result may still fail if absolute regularity is replaced with
the weaker purely nondeterministic assumption. Therefore, the
assumptions in~\cite{Kun71,Kun91,Ste89,OP96} (which implicitly assume
nondegeneracy, though this is not used in the proofs) are genuinely
too weak to yield the desired results.

The results discussed above all assume the classical hidden Markov model
setting illustrated in Figure~\ref{fighmm}(a). Such models are quite
flexible and appear in a wide array of applications~\cite{CMR05}.
Nonetheless, there are many applications in which the need arises for
more general classes of partially observed Markov models. For example,
two common generalizations of the classical hidden Markov model are
illustrated in Figure~\ref{fighmm}(b) and (c). The model of
Figure~\ref{fighmm}(b) is a generalized hidden Markov model~\cite{DS05}
or an autoregressive process with Markov regime~\cite{DMR04}. This
model is similar to a hidden Markov model in that the dynamics of
$(X_n)_{n\ge0}$ do not depend on the observations $(Y_n)_{n\ge0}$;
however, here the observations are not conditionally independent but may
possess their own dynamics. Such models are common in financial
mathematics, where $(Y_n)_{n\ge0}$ might represent a sequence of
investment returns while $(X_n)_{n\ge0}$ models the\vadjust{\goodbreak} state of the
underlying economy. On the other hand, in the model of Figure
\ref{fighmm}(c) there is feedback from the observations to the dynamics
of the unobserved process $(X_n)_{n\ge0}$. Such models arise when the
noise driving the unobserved process and the observation noise are
correlated.

In these more general models, the process $(\Pi_n)_{n\ge0}$ is no
longer Markovian, but the pair $(\Pi_n,Y_n)_{n\ge0}$ is still Markov.
It is therefore natural, and of significant interest for applications,
to investigate the ergodicity of $(\Pi_n,Y_n)_{n\ge0}$ and the
asymptotic stability of $(\Pi_n)_{n\ge0}$ in a more general setting. It
has been shown by Di Masi and Stettner~\cite{DS05} for the model of
Figure~\ref{fighmm}(b), and by Budhiraja~\cite{Bud02} for the model of
Figure~\ref{fighmm}(c), that these problems can be reduced to
establishing the validity of the exchange of intersection and supremum
of $\sigma$-fields along the lines of the earlier approach for classical
hidden Markov models in~\cite{Kun71,Kun91,Ste89,OP96}. The generalization
of the positive results in~\cite{vH10} is far from straightforward,
however.

To illustrate one of the complications that arises in generalized models,
let us consider the setting of Budhiraja~\cite{Bud02}. Budhiraja
considers a model of the form
\[
X_n = f(X_{n-1},Y_{n-1},\xi_{n}),\qquad
Y_n = h(X_n) + \eta_n,
\]
where $(\xi_n)_{n\ge1}$ and $(\eta_n)_{n\ge0}$ are independent
i.i.d.
sequences. It is assumed that $f,h$ are continuous functions and that
$\eta_0$ possesses a bounded and continuous density with respect to some
reference measure $\varphi$. This is evidently a~hidden Markov model
with correlated noise of the type illustrated in Figure~\ref{fighmm}(c).
The main result in~\cite{Bud02} states that if this model admits
a~unique stationary law $\mathbf{P}$ and if $(X_n)_{n\in\mathbb{Z}}$ is
purely nondeterministic, then $(\Pi_n,Y_n)_{n\ge0}$ possesses a unique
invariant measure. Budhiraja's proof contains the same gap as in
\cite{Kun71,Kun91,Ste89}; indeed, the result is clearly erroneous in
light of Example~\ref{exhmm}. Nonetheless, it seems reasonable to
guess that if we assume nondegeneracy of the observations (i.e.,
that the density of $\eta_0$ is strictly positive) and absolute
regularity of the unobserved process, then the result will hold as in
\cite{vH10}. Even this, however, turns out to be false.
%
%ex1.2 #&#
\begin{example}
\label{exchmm}
Define the $\{00,01,10,11\}$-valued process $(X_n)_{n\in\mathbb{Z}}$ and
real-valued process $(Y_n)_{n\in\mathbb{Z}}$ such that
$X_0$ is uniformly distributed in $\{00,01,\allowbreak10,11\}$,
\[
(X_n^1,X_n^2) =
\bigl(X_{n-1}^2,\bigl|X_{n-1}^2-I_{[0,\infty[}(Y_{n-1})\bigr|\bigr),\qquad
Y_n = \eta_n,
\]
where $(\eta_n)_{n\in\mathbb{Z}}$ are i.i.d. $N(0,1)$-distributed
random variables. Then the process
$(X_n,I_{[0,\infty[}(Y_{n-1}))_{n\in\mathbb{Z}}$
has the same law as the classical hidden Markov model of Example
\ref{exhmm}, so stability and unique ergodicity of the filter must fail.
\end{example}

Even though the observations are ostensibly nondegenerate in this
example, the feedback from the observations affects the dynamics of the
unobserved process in a\vadjust{\goodbreak} singular fashion that recreates the problems of
Example~\ref{exhmm}. We thus need at least a different notion of
nondegeneracy in order to rule out such phenomena.

The goal of this paper is to develop a general ergodic and stability
theory for nonlinear filters that subsumes all of the models discussed
above. Indeed, we do not impose any structural assumptions other than
that $(X_n,Y_n)_{n\ge0}$ is a Markov chain that possesses a stationary
law $\mathbf{P}$ [as is illustrated in Figure~\ref{fighmm}(d)]. The main
assumptions of this paper generalize those of~\cite{vH10}; we assume
that the model is:
\begin{longlist}[(2)]
\item[(1)] absolutely regular:
\mbox{$\mathbf{E}(\|\mathbf{P}((X_n,Y_n)\!\in\!\cdv X_0,Y_0)\!-\!
\mathbf{P}((X_n,Y_n)\!\in\!\cdot)\|_{\mathrm{TV}})\!\to\!0$};
\item[(2)] nondegenerate: there exist kernels $P_0,Q$ and
a density
\mbox{$g(x',y',x,y)>0$} so that
$\mathbf{P}((X_{n+1},Y_{n+1})\in A|X_{n},Y_{n})
= \int_A g(X_{n},Y_{n},x,y) P_0(X_{n},dx)
Q(Y_{n},dy)$.
\end{longlist}
The latter assumption states that the dynamics of the observed and
unobserved processes can be made independent (on finite time intervals)
by an equivalent change of measure. It is easily seen that the notion
of nondegenerate observations for the classical hidden Markov model is
a special case of this assumption; on the other hand, the present
assumption also rules out the phenomenon observed in
Example~\ref{exchmm}. This general nondegeneracy property appears to be
precisely the right assumption required to generalize the results of
\cite{vH10}, and seems very natural in view of Examples~\ref{exhmm} and
\ref{exchmm}. The absolute regularity assumption on
$(X_n,Y_n)_{n\in\mathbb{Z}}$ can in fact be weakened somewhat; see
Sections~\ref{secmainres} and~\ref{secsuffcond} for a precise
statement.

With the above assumptions in place, we will show that Kunita's exchange
of intersection and supremum of $\sigma$-fields is permitted in our
setting, and we can consequently develop general asymptotic stability
and unique ergodicity results. The intuition behind the proofs is
similar in spirit to the classical hidden Markov model setting in
\cite{vH10,vH11}, and we refer to those papers for a~discussion of the
basic ideas. Nonetheless, to our surprise, key parts of the proofs in
\cite{vH10} break down completely in the generalized setting of this
paper and almost all arguments in~\cite{vH10} require substantial
modification, as we can no longer exploit many simplifying properties
that hold trivially in classical hidden Markov models. The proofs in
the present paper rely on the ergodic properties of nondegenerate Markov
chains that are developed in Section~\ref{secnondeg} below. Though this
paper is almost entirely self-contained, the reader may find it helpful
to familiarize herself first with the simpler setting of~\cite{vH10}.

This paper is organized as follows. Section~\ref{secmain} introduces the
general model used throughout the paper and states our main results. We
also give useful sufficient conditions for the models in Figure
\ref{fighmm}(a)--(d). Section~\ref{secnondeg} develops the
ergodic properties of nondegenerate Markov chains that play a central role
in our proofs. Sections~\ref{secexchg}--\ref{seclambdainv} are devoted
to the proofs of our main results. Appendices~\ref{secapp} and
\ref{secnotationlist} collect auxiliary results and a notation list
that is used throughout the paper.

%s2 #&#
\section{Preliminaries and main results}
\label{secmain}

%s2.1 #&#
\subsection{The canonical setup}

Throughout this paper we consider the bivariate stochastic process
$(X_n,Y_n)_{n\in\mathbb{Z}}$, where $X_n$ takes values in the Polish
space~$E$ and $Y_n$ takes values in the Polish space $F$.\label{page7}
We realize this process on the canonical path space
$\Omega=\Omega^X\times\Omega^Y$ with
$\Omega^X= E^{\mathbb{Z}}$ and $\Omega^Y=F^{\mathbb{Z}}$, such that
$X_n(x,y)=x(n)$ and $Y_n(x,y)=y(n)$.\label{page7bb}
%component}%
%component}%
Denote by $\mathcal{F}$ the Borel $\sigma$-field of~$\Omega$, and
define\label{pageBorelsigma}
\[
\mathcal{F}_I^X=\sigma\{X_k\dvtx k\in I\},\qquad
\mathcal{F}_I^Y=\sigma\{Y_k\dvtx k\in I\},\qquad
\mathcal{F}_I=\mathcal{F}_I^X\vee\mathcal{F}_I^Y
\]
for $I\subset\mathbb{Z}$.\label{pageFF}
For simplicity of notation, we define the natural filtrations
\[
\mathcal{F}_n^X=\mathcal{F}_{]-\infty,n]}^X,\qquad
\mathcal{F}_n^Y=\mathcal{F}_{]-\infty,n]}^Y,\qquad
\mathcal{F}_n=\mathcal{F}_{]-\infty,n]} \qquad(n\in\mathbb{Z})
\]
and the $\sigma$-fields\label{pageFFY}
\[
\mathcal{F}^X=\mathcal{F}^X_{\mathbb{Z}},\qquad
\mathcal{F}^Y=\mathcal{F}^Y_{\mathbb{Z}},\qquad
\mathcal{F}^X_+=\mathcal{F}^X_{[0,\infty[},\qquad
\mathcal{F}^Y_+=\mathcal{F}^Y_{[0,\infty[}.
\]
Finally, we denote by $Y$ the $F^{\mathbb{Z}}$-valued random variable
$(Y_k)_{k\in\mathbb{Z}}$, and the canonical shift $\Theta\dvtx\Omega\to
\Omega$
is defined as $\Theta(x,y)(m) = (x(m+1),y(m+1))$.\label{page7b}

For any Polish space $Z$, we denote by $\mathcal{B}(Z)$ its Borel
$\sigma$-field and by $\mathcal{P}(Z)$ the space of all probability
measures on $Z$ endowed with the weak convergence topology [thus
$\mathcal{P}(Z)$ is again Polish].\label{pageBorel}
Let us recall that any probability kernel
$\rho\dvtx Z\times\mathcal{B}(Z')\to[0,1]$ may be equivalently viewed as
a $\mathcal{P}(Z')$-valued random variable $z\mapsto\rho(z, \cdot)$ on
$(Z,\mathcal{B}(Z))$. For notational convenience, we will implicitly
identify probability kernels and random probability measures in the
sequel.

%s2.2 #&#
\subsection{The model}

The basic model of this paper is defined by a Markov transition kernel
$P\dvtx E\times F\times\mathcal{B}(E\times F)\to[0,1]$ and a $P$-invariant
probability measure $\pi$ on $(E\times F, \mathcal{B}(E\times F))$, which
we presume to be fixed throughout the paper. We now define the
probability measure $\mathbf{P}$ on $(\Omega,\mathcal{F})$ such that,
under $\mathbf{P}$, the process $(X_n,Y_n)_{n\in\mathbb{Z}}$ is the
stationary Markov chain with transition kernel $P$ and stationary
distribution $\pi$.\label{page7a}
We interpret $Y_n$ to be the \textit{observable
component} of the model, while $X_n$ is the \textit{unobservable
component}.\label{page7aaa}

As $(X_n,Y_n)_{n\in\mathbb{Z}}$ is a stationary Markov chain under
$\mathbf{P}$, the reverse time process $(X_{-n},Y_{-n})_{n\in\mathbb{Z}}$
is again a stationary Markov chain. We fix throughout the paper a version
$P'\dvtx E\times F\times\mathcal{B}(E\times F)\to[0,1]$ of the regular
conditional probability $\mathbf{P}((X_{-1},Y_{-1})\in\cdv X_0,Y_0)$.
Thus, by construction, the process $(X_{-n},Y_{-n})_{n\in\mathbb{Z}}$
is a
stationary Markov chain with transition kernel~$P'$ and invariant measure
$\pi$.\label{page7aa}
%$(X_{-n},Y_{-n})_{n\in\mathbb{Z}}$}%

In addition to the probability measure $\mathbf{P}$, we introduce
the probability kernel $\mathbf{P}^\cdot\dvtx E\times F
\times\mathcal{F}\to[0,1]$ with the following
properties:\label{pageCond}
under $\mathbf{P}^{z,w}$,
\begin{longlist}[(2)]
\item[(1)] $(X_n,Y_n)_{n\ge0}$ is Markov with transition kernel
$P$ and initial measure $\delta_z\otimes\delta_w$;\vadjust{\goodbreak}
\item[(2)] $(X_{-n},Y_{-n})_{n\ge0}$ is Markov with transition kernel
$P'$ and initial measure $\delta_z\otimes\delta_w$;
\item[(3)] $(X_n,Y_n)_{n\ge0}$ and $(X_{-n},Y_{-n})_{n\ge0}$ are independent.
\end{longlist}
Clearly $\mathbf{P}^{z,w}$ is a version of the regular conditional
probability $\mathbf{P}( \cdv X_0,Y_0)$. Finally, for any probability
measure $\nu$ on $(E\times F, \mathcal{B}(E\times F))$, we define
\[
\mathbf{P}^\nu(A) = \int I_A(x,y) \mathbf{P}^{z,w}(dx,dy) \nu(dz,dw)
\qquad\mbox{for all }A\in\mathcal{F}.
\]
Note, in particular, that $\mathbf{P}^\pi$ coincides with $\mathbf{P}$
by construction.\label{pagePP}

%s2.3 #&#
\subsection{The nonlinear filter}
\label{secintrofilt}

As $X_n$ is not directly observable, we are interested in the conditional
distribution of $X_n$ given the history of observations to date
$Y_0,\ldots,Y_n$. To this end, we define for every probability measure
$\mu$ on $E\times F$ and $n\ge0$ the \textit{nonlinear filter}
$\Pi_n^\mu\dvtx\Omega^Y\times\mathcal{B}(E)\to[0,1]$ to be a version of the
regular conditional probability
$\mathbf{P}^\mu(X_n\in\cdv \mathcal{F}^Y_{[0,n]})$.\label{pagea8}
The nonlinear filter is the central object of interest throughout this
paper.

We now state some basic properties of the nonlinear filter. The first
property establishes that the filter can be computed recursively.
%
%le2.1 #&#
\begin{lem}
\label{lemintrofilt}
There is a measurable map $U\dvtx\mathcal{P}(E)\times F\times F\to
\mathcal{P}(E)$ such that $\Pi_n^\mu= U(\Pi_{n-1}^\mu,Y_{n-1},Y_n)$
$\mathbf{P}^\mu$-a.s. for every $n\ge1$ and $\mu\in\mathcal
{P}(E\times F)$.\label{page8}
%section~\ref{secpairm})}%
\end{lem}
%
%re2.2 #&#
\begin{rem}
In the proof of our main results, it will be convenient to assume that the
identity $\Pi_n^\mu= U(\Pi_{n-1}^\mu,Y_{n-1},Y_n)$ holds everywhere on
$\Omega^Y$ and not just $\mathbf{P}^\mu$-a.s. This corresponds to the
choice of a particular version of the nonlinear filter. However, as none
of our results will depend on the choice of version of the filter,
there is clearly no loss of generality in fixing such a convenient version
for the purposes of our proofs, as we will do in Section~\ref{secfiltstab}.
\end{rem}

We now consider $(\Pi_n^\mu)_{n\ge0}$ as a $\mathcal{P}(E)$-valued
stochastic process. The second property establishes that this
measure-valued process inherits certain Markovian properties from the
underlying model $(X_n,Y_n)_{n\ge0}$.
%
%le2.3 #&#
\begin{lem}
\label{lemintromarkov}
There exist Markov transition kernels $\Gamma$ on $\mathcal{P}(E)\times F$
and~$\Lambda$ on $\mathcal{P}(E)\times E\times F$ such that
the following hold:\label{page8a}
%0}$ (Lemma~\ref{lemintromarkov}, cf. section~\ref{secpairm})}%
for every $\mu\in\mathcal{P}(E\times F)$,
\begin{longlist}[(2)]
\item[(1)] $(\Pi_n^\mu,Y_n)_{n\ge0}$ is a Markov chain under
$\mathbf{P}^\mu$ with transition kernel $\Gamma$; and
\item[(2)] $(\Pi_n^\mu,X_n,Y_n)_{n\ge0}$ is a Markov chain under
$\mathbf{P}^\mu$ with transition kernel~$\Lambda$.
\end{longlist}
\end{lem}

For any $\mathsf{m}\in\mathcal{P}(\mathcal{P}(E)\times F)$, define
the \textit{barycenter} $b\mathsf{m}\in\mathcal{P}(E\times F)$
as\label{pageBary}
\[
b\mathsf{m}(A\times B) =
\int\nu(A) I_B(w) \mathsf{m}(d\nu,dw).
\]
We finally state some properties of $\Gamma$- and $\Lambda$-invariant
measures.\vadjust{\goodbreak}
%
%le2.4 #&#
\begin{lem}
\label{lemintrobary}
For any $\Gamma$-invariant probability measure
$\mathsf{m}\in\mathcal{P}(\mathcal{P}(E)\times F)$, the barycenter
$b\mathsf{m}$ is a $P$-invariant probability measure. Conversely,
there exists at least one $\Gamma$-invariant probability measure with
barycenter $\pi$.

Similarly, for any $\Lambda$-invariant probability measure
$\mathsf{M}\in\mathcal{P}(\mathcal{P}(E)\times E\times F)$, the
marginal $\mathsf{M}(\mathcal{P}(E)\times\cdot)$ is a $P$-invariant
probability measure. Conversely, there exists at least one
$\Lambda$-invariant probability measure with marginal $\pi$.
\end{lem}

In general, there may be multiple $\Gamma$-invariant measures with
barycenter $\pi$, etc. Our main results will establish uniqueness
under suitable assumptions.
%
%re2.5 #&#
\begin{rem}
For the purposes of this paper it suffices to establish the above
results for the case where Assumption~\ref{asptnondeg} below is assumed
to hold. In this setting, these results will be proved in Sections
\ref{secpairm} and~\ref{sectriplem}. In fact, the results in this
subsection hold very generally as stated without any further
assumptions, but the proofs in the general setting are somewhat more
abstract. Such generality will not be needed in this paper, and we
therefore leave the generalization of the proofs (along the lines of
\cite{vH11}, Appendix A.1) to the interested reader.
\end{rem}

%s2.4 #&#
\subsection{Main results}
\label{secmainres}

We begin by introducing the fundamental model assumptions that are
required by our main results. Let us emphasize that we will at no point
in the paper automatically assume that any of these assumptions is in
force; all assumptions will be imposed explicitly where they are needed.
Some useful sufficient conditions will be given in Section~\ref{secsuffcond}
below.

%as2.6 #&#
\begin{aspt}[(Marginal ergodicity)]
\label{asptergodic} The following holds:
\[
\int
\|\mathbf{P}^{z,w}(X_n\in\cdot)-\mathbf{P}(X_n\in\cdot)\|_{\mathrm{TV}}
\pi(dz,dw)
\xrightarrow{n\to\infty}0.
\]
\end{aspt}
%
%as2.7 #&#
\begin{aspt}[(Reversed marginal ergodicity)]
\label{asptergodicrv} The following holds:
\[
\int
\|\mathbf{P}^{z,w}(X_{-n}\in\cdot)-\mathbf{P}(X_{-n}\in\cdot)\|
_{\mathrm{TV}}
\pi(dz,dw)
\xrightarrow{n\to\infty}0.
\]
\end{aspt}
%
%as2.8 #&#
\begin{aspt}[(Nondegeneracy)]
\label{asptnondeg}
There exist transition probability kernels
$P_0\dvtx E\times\mathcal{B}(E)\to[0,1]$ and $Q\dvtx F\times\mathcal{B}(F)\to[0,1]$
such that
\[
P(z,w,dz',dw')=g(z,w,z',w')P_0(z,dz')Q(w,dw')
\]
for some strictly positive measurable function
$g\dvtx E\times F\times E\times F\to\ ]0,\infty[$.\label{page9}
%Q$ (Assumption~\ref{asptnondeg})}%
\end{aspt}

We now proceed to state the main results of this paper. Our results
address in turn each of the problems discussed in the
\hyperref[intro]{Introduction}:
the exchange of intersection and supremum of $\sigma$-fields,
asymptotic stability of the nonlinear filter and unique ergodicity
of the processes $(\Pi_n^\mu,Y_n)_{n\ge0}$ and
$(\Pi_n^\mu,X_n,Y_n)_{n\ge0}$.\vadjust{\goodbreak}

Our first result establishes the validity of Kunita's exchange of
intersection and supremum, and its time-reversed cousin, in the
generalized setting of this paper.
%
%th2.9 #&#
\begin{theorem}
\label{thmexchg} Suppose that Assumptions
\ref{asptergodic}--\ref{asptnondeg} are in force. Then
\[
\bigcap_{n\ge0} \mathcal{F}^Y_+\vee\mathcal{F}^X_{[n,\infty[}=
\mathcal{F}^Y_+ \quad\mbox{and}\quad
\bigcap_{n\ge0} \mathcal{F}^Y_0\vee\mathcal{F}^X_{-n}=
\mathcal{F}^Y_0 ,\qquad \mathbf{P}\mbox{-a.s.}
\]
\end{theorem}

Our second result concerns filter stability which can be established
in our setting (as in~\cite{vH10}) in a very strong sense: pathwise
and in the total variation topology.
%
%th2.10 #&#
\begin{theorem}
\label{thmfilter} Suppose that Assumptions
\ref{asptergodic}--\ref{asptnondeg} are in force. Let
$\mu$ be a probability measure on $E\times F$ such that
$\mu(E\times\cdot)\ll\pi(E\times\cdot)$ and
\[
\mathbf{E}^\mu\bigl(
\bigl\|\mathbf{P}^\mu(X_n\in\cdv Y_0)-\mathbf{P}(X_n\in\cdot)
\bigr\|_{\mathrm{TV}}\bigr)
\xrightarrow{n\to\infty}0.
\]
Then $\|\Pi_n^\mu-\Pi_n^\pi\|_{\mathrm{TV}}\xrightarrow{n\to\infty}0$
$\mathbf{P}^\mu$-a.s. [and $\mathbf{P}$-a.s. if
$\mu(E\times\cdot)\sim\pi(E\times\cdot)$].
\end{theorem}
%
%re2.11 #&#
\begin{rem}
The assumptions of Theorem~\ref{thmfilter} may be more intuitive when
phrased in terms of the filtering recursion in Lemma~\ref{lemintrofilt}.
Let $\rho\dvtx F\times\mathcal{B}(E)\to[0,1]$ be a probability kernel,
and define the random measures $(\Pi_n)_{n\ge0}$ by the recursion
\[
\Pi_0 = \rho(Y_0, \cdot),\qquad
\Pi_n = U(\Pi_{n-1},Y_{n-1},Y_n).
\]
Suppose that the dynamics of $(X_n)_{n\ge0}$ are such that the random
initial law~$\rho$ is in the domain of attraction of the stationary
distribution $\pi$ in the sense that
\[
\bigl\|\mathbf{P}^{\rho(w, \cdot)\otimes\delta_w}(X_n\in\cdot)-
\mathbf{P}(X_n\in\cdot)\bigr\|_{\mathrm{TV}}\xrightarrow{n\to\infty}0
\qquad\mbox{in }\pi(E\times dw)\mbox{-probability}.
\]
Then $\|\Pi_n-\Pi_n^\pi\|_{\mathrm{TV}}\xrightarrow{n\to\infty}0$
$\mathbf{P}$-a.s. Indeed, this follows immediately from Theorem
\ref{thmfilter} by setting $\mu(dz,dw) = \rho(w,dz)\pi(E\times dw)$.
Therefore, we may interpret Theorem~\ref{thmfilter} as follows:
the filtering recursion of Lemma~\ref{lemintrofilt} is asymptotically
stable inside the domain of attraction of the stationary distribution.

The result of Theorem~\ref{thmfilter} is easily extended to show
\mbox{$\|\Pi_n^\mu-\Pi_n^\nu\|_{\mathrm{TV}}\xrightarrow{n\to\infty}0$}
$\mathbf{P}^\gamma$-a.s. whenever all three initial
measures $\mu,\nu,\gamma$ are in the domain of attraction of the
stationary distribution in the above sense, using
Corollary~\ref{coryinfhabscont} below.
\end{rem}

Our third result concerns uniqueness of the $\Gamma$-invariant measure.
%
%th2.12 #&#
\begin{theorem}
\label{thminvmeas}
Suppose that Assumptions~\ref{asptergodic}--\ref{asptnondeg} are in
force. Then there exists a unique $\Gamma$-invariant probability measure
with barycenter $\pi$. In particular, if $P$ has a unique invariant
probability measure, then so does $\Gamma$.\vadjust{\goodbreak}
\end{theorem}

Our fourth result concerns uniqueness of the $\Lambda$-invariant
measure. The situation here is a little more complicated; Assumptions
\ref{asptergodic}--\ref{asptnondeg} only ensure uniqueness within a
restricted class of measures (cf.~\cite{Kun91}), while a somewhat
stronger variant of Assumption~\ref{asptergodic} yields uniqueness in
the class of all probability measures.
%
%th2.13 #&#
\begin{theorem}
\label{thmlambdainv}
Suppose that Assumptions~\ref{asptergodic}--\ref{asptnondeg} hold.
Then there exists a unique $\Lambda$-invariant probability measure with
marginal $\pi$ on $E\times F$ in the class
\begin{eqnarray*}
&&\biggl\{\mathsf{M}\in\mathcal{P}\bigl(\mathcal{P}(E)\times E\times F\bigr)\mbox{:
for every }
A\in\mathcal{B}(\mathcal{P}(E)),
B\in\mathcal{B}(E),
C\in\mathcal{B}(F), \\
&&\hspace*{90pt}\mathsf{M}(A\times B\times C) =
\int\nu(B) I_{A\times C}(\nu,w) \mathsf{M}(d\nu,dz,dw)
\biggr\}.
\end{eqnarray*}
If, in addition, we have
\[
\int
\|\mathbf{P}^{z,w}(X_n\in\cdot)-\mathbf{P}(X_n\in\cdot)
\|_{\mathrm{TV}} \mu(dz,dw)\xrightarrow{n\to\infty}0
\]
for every probability measure $\mu$ on $E\times F$ such that
$\mu(E\times\cdot)=\pi(E\times\cdot)$,
then there exists a unique $\Lambda$-invariant probability measure with
marginal~$\pi$ among all probability measures in
$\mathcal{P}(\mathcal{P}(E)\times E\times F)$.
If we assume even further that~$P$ has a unique invariant probability
measure, then so does $\Lambda$.
\end{theorem}

The following sections are devoted to the proofs of these results:
Theorems~\ref{thmexchg},~\ref{thmfilter},~\ref{thminvmeas} and
\ref{thmlambdainv} are proved in Sections~\ref{secexchg},
\ref{secfiltstab},~\ref{secinvmeas} and~\ref{seclambdainv},
respectively.

%s2.5 #&#
\subsection{Sufficient conditions}
\label{secsuffcond}

Our main results rely on the fundamental Assumptions
\ref{asptergodic}--\ref{asptnondeg}. In most applications, the form
of the transition kernel~$P$ is explicitly (or semi-explicitly) given.
Existence and uniqueness of an invariant measure $\pi$ and the
ergodicity Assumption~\ref{asptergodic} can often be verified in terms
of $P$ only (cf.~\cite{MT09}), while the nondegeneracy Assumption
\ref{asptnondeg} can be read off directly from the explicit form of
$P$. On the other hand, explicit expressions for the invariant measure
$\pi$ or the reversed transition kernel $P'$ are often not available, so
that Assumption~\ref{asptergodicrv} may be difficult to verify
directly. The goal of this section is to provide sufficient conditions
for our main results that are easily verified in practice.

%s2.5.1 #&#
\subsubsection{General sufficient conditions}

Our main sufficient condition is absolute regularity (cf.~\cite{VR59}),
of the process $(X_n,Y_n)_{n\in\mathbb{Z}}$, which was the assumption
stated in the \hyperref[intro]{Introduction}. This is slightly stronger than
Assumptions~\ref{asptergodic} and~\ref{asptergodicrv}, but has the benefit that it
is automatically time-reversible and therefore easily verifiable.
%
%le2.14 #&#
\begin{lem}
\label{lemabsreg}
Suppose that $(X_n,Y_n)_{n\in\mathbb{Z}}$ is absolutely regular,
\[
\int\bigl\|\mathbf{P}^{z,w}\bigl((X_n,Y_n)\in\cdot\bigr)
-\pi\bigr\|_{\mathrm{TV}} \pi(dz,dw)
\xrightarrow{n\to\infty}0.
\]
Then both Assumptions~\ref{asptergodic} and
\ref{asptergodicrv} hold true.
\end{lem}
\begin{pf}
Absolute regularity trivially yields Assumption~\ref{asptergodic}.
On the other hand, the absolute regularity property of a stationary
Markov chain is invariant under time reversal by~\cite{vH10},
Proposition 4.4, so that Assumption~\ref{asptergodicrv} follows.
\end{pf}

Similarly, the convergence assumption in Theorem~\ref{thmfilter} also
admits a~slight\-ly stronger but potentially more easily verified
counterpart.
%
%le2.15 #&#
\begin{lem}
\label{lemabsregcvg}
Suppose that Assumption~\ref{asptergodic} holds. Let
$\mu$ be a probability measure on $E\times F$ such that
$\|\mathbf{P}^\mu((X_n,Y_n)\in\cdot)-\pi\|_{\mathrm{TV}}\to0$
as $n\to\infty$. Then
\[
\mathbf{E}^\mu\bigl(
\bigl\|\mathbf{P}^\mu(X_n\in\cdv Y_0)-\mathbf{P}(X_n\in\cdot)
\bigr\|_{\mathrm{TV}}\bigr)\xrightarrow{n\to\infty}0.
\]
\end{lem}
\begin{pf}
Define the quantity
\[
\Delta_k(x,y) =
\|\mathbf{P}^{x,y}(X_{k}\in\cdot)-\mathbf{P}(X_{k}\in\cdot)
\|_{\mathrm{TV}}.
\]
By the stationarity of $\mathbf{P}$, the Markov property
and $\|\Delta_k-1\|_\infty\le1$, we can estimate
\begin{eqnarray*}
&&\mathbf{E}^\mu\bigl(
\bigl\|\mathbf{P}^\mu(X_{n+k}\in\cdv Y_0)
-\mathbf{P}(X_{n+k}\in\cdot)\bigr\|_{\mathrm{TV}}\bigr) \\
&&\qquad \le
\mathbf{E}^\mu\bigl(
\|\mathbf{P}^{X_n,Y_n}(X_{k}\in\cdot)-\mathbf{P}(X_{k}\in\cdot)
\|_{\mathrm{TV}}\bigr) \\
&&\qquad =
\mathbf{E}(\Delta_k(X_n,Y_n)) +
\bigl\{\mathbf{E}^\mu\bigl(\Delta_k(X_n,Y_n)-1\bigr)-
\mathbf{E}\bigl(\Delta_k(X_n,Y_n)-1\bigr)\bigr\} \\
&&\qquad \le
\mathbf{E}\bigl(
\|\mathbf{P}^{X_0,Y_0}(X_{k}\in\cdot)-\mathbf{P}(X_{k}\in\cdot)
\|_{\mathrm{TV}}\bigr)+
\bigl\|\mathbf{P}^\mu\bigl((X_n,Y_n)\in\cdot\bigr)-\pi\bigr\|_{\mathrm{TV}}.
\end{eqnarray*}
This expression converges to zero as $k,n\to\infty$ by our assumptions.
\end{pf}

%s2.5.2 #&#
\subsubsection{Generalized hidden Markov models}

We now consider the special case where the underlying model
$(X_n,Y_n)_{n\in\mathbb{Z}}$ is a generalized hidden Markov model,
whose dependence structure is illustrated in Figure~\ref{fighmm}(b).
Under Assumption~\ref{asptnondeg}, this dependence structure is
enforced by the additional requirement that
\[
\int g(z,w,z',w') Q(w,dw') = 1 \qquad\mbox{for all }w\in F,
z,z'\in E.
\]
This implies that $(X_n)_{n\in\mathbb{Z}}$ is itself Markovian under
$\mathbf{P}$ with transition kernel~$P_0$, and the probability measure
$\pi_0=\pi( \cdot\times F)$ must then be $P_0$-invariant. In this
setting, it suffices to consider the ergodic properties of the
unobserved process, provided that the reference transition kernel
$Q(w,dw')$ does not depend on $w$.
%
%le2.16 #&#
\begin{lem}
\label{lemgenhmm}
$\!\!\!$Suppose that Assumption~\ref{asptnondeg} holds with
$Q(w,dw')= \varphi(dw')$ for some probability measure $\varphi$ on $F$,
and that $(X_n,Y_n)_{n\in\mathbb{Z}}$ is a generalized hidden Markov model
in the above sense. If $(X_n)_{n\in\mathbb{Z}}$ is absolutely regular
\[
\int\| P_0^n(z, \cdot)-\pi_0\|_{\mathrm{TV}}
\pi_0(dz)\xrightarrow{n\to\infty}0,
\]
then both Assumptions~\ref{asptergodic} and
\ref{asptergodicrv} hold true.
\end{lem}
\begin{pf}
We reduce to the case of Lemma~\ref{lemabsreg}. A stationary Markov
chain is absolutely regular if and only if for almost every pair of
initial conditions, there is a finite time $n$ at which the laws of the
chain are not mutually singular (e.g., this is a special case of
Theorem~\ref{thmmcre} below). Therefore, our assumption implies that
for $\pi_0\otimes\pi_0$-a.e. $(z,z')$, there is an $n\ge0$ such that
$P_0^n(z, \cdot)$ and $P_0^n(z', \cdot)$ are not mutually
singular. But as $Q(w,dw')=\varphi(dw')$ and by Assumption
\ref{asptnondeg}, we have $P^n(z,w, \cdot)\sim P_0^n(z, \cdot)
\otimes\varphi$ and $P^n(z',w', \cdot)\sim P_0^n(z', \cdot)
\otimes\varphi$ for every $z,w,z',w'$. It follows that for
$\pi\otimes\pi$-a.e. $((z,w),(z',w'))$ there is an $n\ge0$ such that
$P^n(z,w, \cdot)$ and $P^n(z',w', \cdot)$ are not mutually singular.
We have therefore shown that the absolutely regularity assumption of Lemma
\ref{lemabsreg} holds.
\end{pf}
%
%re2.17 #&#
\begin{rem}
By the generalized hidden Markov structure
$\mathbf{P}^{z,w}(X_n\in\cdot) =P_0^n(z, \cdot)$ is independent
of $w$, so that Assumption~\ref{asptergodic} follows immediately from
the absolute regularity of $(X_n)_{n\in\mathbb{Z}}$. Unfortunately, the
generalized hidden Markov property is not invariant under time reversal,
so this argument does not guarantee that Assumption~\ref{asptergodicrv}
holds. The additional assumption that $Q(w,dw')=\varphi(dw')$ allows us to
circumvent this problem by reducing to the case of Lemma~\ref{lemabsreg}.
\end{rem}

We also have a counterpart of Lemma~\ref{lemabsregcvg}.
%
%le2.18 #&#
\begin{lem}
\label{lemghmmcvg}
Suppose the assumptions of Lemma~\ref{lemgenhmm} hold. Let
$\mu$ be a~probability measure on $E\times F$ so that
$\|\mu( \cdot\times F)P_0^n-\pi_0\|_{\mathrm{TV}}\to0$ as $n\to\infty$. Then
\[
\mathbf{E}^\mu\bigl(
\bigl\|\mathbf{P}^\mu(X_n\in\cdv Y_0)-\mathbf{P}(X_n\in\cdot)
\bigr\|_{\mathrm{TV}}\bigr)\xrightarrow{n\to\infty}0.
\]
\end{lem}
\begin{pf}
We reduce to the case of Lemma~\ref{lemabsregcvg}. As
$Q(w,dw')=\varphi(dw')$, we obtain $\pi\sim\pi_0\otimes\varphi$ and
$\mu P^n \sim\mu( \cdot\times F)P_0^n\otimes\varphi$ for all $n>0$
by Assumption~\ref{asptnondeg}. Choose $S_n\in\mathcal{B}(E)$ such
that $\mu( \cdot\times F)P_0^n( {\cdot}\cap S_n)\ll\pi_0$ and
$\pi_0(S_n^c)=0$ for all $n$ [so~$S_n$ defines the Lebesgue decomposition
of $\mu( \cdot\times F)P_0^n$ with respect to $\pi_0$]. Then
clearly
$\mu P^n( \cdot\cap S_n\times F)\ll\pi$ and $\pi(S_n^c\times F)=0$.
Therefore
\begin{eqnarray*}
\|\mu P^{k+n}-\pi\|_{\mathrm{TV}} &\le&
\mu P^n(S_n\times F) \|\nu_n P^k-\pi\|_{\mathrm{TV}} +
2 \mu P^n(S_n^c\times F) \\
&\le&\|\nu_n P^k-\pi\|_{\mathrm{TV}} +
\|\mu( \cdot\times F)P_0^n-\pi_0\|_{\mathrm{TV}},\vadjust{\goodbreak}
\end{eqnarray*}
where we have defined $\nu_n= \mu P^n( \cdot\cap
S_n\times F)/\mu P^n(S_n\times F)$. But as $(X_n,Y_n)_{n\in\mathbb{Z}}$
is absolutely regular (cf. Lemma~\ref{lemgenhmm}) and
$\nu_n\ll\pi$, the first term converges to zero as $k\to\infty$.
Letting $n\to\infty$ and applying
Lemma~\ref{lemabsregcvg} yields the result.~%
\end{pf}

%s2.5.3 #&#
\subsubsection{Hidden Markov models with correlated noise}

We now turn to the special case where the underlying model
$(X_n,Y_n)_{n\in\mathbb{Z}}$ is a hidden Markov model with correlated
noise, whose dependence structure is illustrated in Figure
\ref{fighmm}(c). Under Assumption~\ref{asptnondeg}, this
dependence structure is enforced by the following requirement:
there is a probability measure $\varphi$ on $F$ such that
$Q(w,dw')=\varphi(dw')$, and there are measurable functions
\mbox{$g_X\dvtx E\times F\times E\to\mathbb{R}_+$} and
$g_Y\dvtx E\times F\to\mathbb{R}_+$ such that
\[
g(z,w,z',w') = g_X(z,w,z') g_Y(z',w'),\qquad
\int g_Y(z,w) \varphi(dw) = 1.
\]
Unlike in the case of a generalized hidden Markov model, in the present
model the probabilities $\mathbf{P}^{z,w}(X_n\in\cdot)$ do depend
on $w$. Nonetheless, in the present case the unobserved process
$(X_n)_{n\in\mathbb{Z}}$ is still Markov under the stationary measure
$\mathbf{P}$ with respect to its own filtration, with transition
kernel~$\tilde P_0$ given for $A\in\mathcal{B}(E)$ by
\[
\tilde P_0(z,A) = \int P(z,w,A\times F) g_Y(z,w) \varphi(dw).
\]
To see\vspace*{1pt} this, note that
$\pi(dz,dw) = g_Y(z,w) \pi(dz\times F) \varphi(dw)$ by our assumption
on $P$ and $\pi P =\pi$, so we can compute
$\mathbf{P}(X_{n+1}\in A|\mathcal{F}^X_n) = \tilde P_0(X_n,A)$.

%re2.19 #&#
\begin{rem}
Unlike in the case of a generalized hidden Markov model, where
$Q(w,dw')=\varphi(dw')$ is an additional assumption, in the
present setting the assumption $Q(w,dw')=\varphi(dw')$ entails
no loss of generality. Indeed, the hidden Markov structure
with correlated noise can be generally formulated by the requirement
that $P(z,w,dz',dw') = P_X(z,w,dz') P_Y(z',dw')$ for some probability
kernels $P_X$ and $P_Y$. It is easily seen that any such model
that also satisfies Assumption~\ref{asptnondeg} must have the above
form for a suitable choice of $\varphi$.
\end{rem}

The idea is now that in the present setting, it suffices to consider
the ergodic properties of the unobserved process (i.e., the
transition kernel $\tilde P_0$).

%le2.20 #&#
\begin{lem}
\label{lemchmm}
Suppose that Assumption~\ref{asptnondeg} holds
and that $(X_n,Y_n)_{n\in\mathbb{Z}}$ is a hidden Markov model with
correlated noise in the above sense. If also
\[
\int\| \tilde P_0^n(z, \cdot)-\pi_0\|_{\mathrm{TV}}
\pi_0(dz)\xrightarrow{n\to\infty}0,
\]
where $\pi_0=\pi( \cdot\times F)$,
then both Assumptions~\ref{asptergodic} and
\ref{asptergodicrv} hold true.
\end{lem}
\begin{pf}
Note that for all $(z,w)\in E\times F$, $B\in\mathcal{B}(E\times F)$ and
$n\ge1$ we have
\[
P^n(z,w,B) =
\int\biggl\{\int_B
\tilde P_0^{n-1}(z',d\tilde z) g_Y(\tilde z,\tilde w)
\varphi(d\tilde w)
\biggr\}
g_X(z,w,z') P_0(z,dz').
\]
Therefore, we have for $n\ge1$
\[
\int\|P^n(z,w, \cdot)
-\pi\|_{\mathrm{TV}} \pi(dz,dw) \le
\int\| \tilde P_0^{n-1}(z, \cdot)-\pi_0\|_{\mathrm{TV}}.
\]
The result now follows directly from Lemma~\ref{lemabsreg}.
\end{pf}

In the present setting (as in the case of a classical hidden Markov
model), the most natural initial measures $\mu$ are those that are
compatible with the observation model in the sense that $\mu(dz,dw) =
g_Y(z,w) \mu_0(dz) \varphi(dw)$ for some probability measure $\mu_0$
on $E$. This yields the following counterpart of Lemma
\ref{lemghmmcvg}, whose proof (by reduction to Lemma
\ref{lemabsregcvg}) is trivial and is therefore omitted.
%
%le2.21 #&#
\begin{lem}
Suppose the assumptions of Lemma~\ref{lemchmm} hold. Let
$\mu_0$ be a probability measure on $E$ such that
$\|\mu_0\tilde P_0^n-\pi_0\|_{\mathrm{TV}}\to0$ as $n\to\infty$. Then
\[
\mathbf{E}^\mu\bigl(
\bigl\|\mathbf{P}^\mu(X_n\in\cdv Y_0)-\mathbf{P}(X_n\in\cdot)
\bigr\|_{\mathrm{TV}}\bigr)\xrightarrow{n\to\infty}0,
\]
where we have defined $\mu(dz,dw) =
g_Y(z,w) \mu_0(dz) \varphi(dw)$.
\end{lem}
%
%re2.22 #&#
\begin{rem}
Let us note that in all of the special cases discussed above the process
$(X_n,Y_n)_{n\in\mathbb{Z}}$ is absolutely regular so that Assumptions
\ref{asptergodic} and~\ref{asptergodicrv} hold by virtue of Lemma
\ref{lemabsreg}. Absolute regularity of $(X_n,Y_n)_{n\in\mathbb{Z}}$
is not necessary, however, for Assumptions~\ref{asptergodic} and
\ref{asptergodicrv} to hold. For example, in the trivial case that
Assumption~\ref{asptnondeg} holds with $g\equiv1$, it is easily seen
that Assumptions~\ref{asptergodic} and~\ref{asptergodicrv} hold if and
only if the unobserved process $(X_n)_{n\in\mathbb{Z}}$ is absolutely
regular, while the pair process $(X_n,Y_n)_{n\in\mathbb{Z}}$ need not
even be ergodic [e.g., when $Q(w,dw')=\delta_{w}(dw')$]. Thus
Assumptions~\ref{asptergodic} and~\ref{asptergodicrv} are strictly
weaker than the absolute regularity of the pair process
$(X_n,Y_n)_{n\in\mathbb{Z}}$. Nonetheless, the latter assumption is very
mild and will likely hold in most applications of practical interest.
\end{rem}

%s3 #&#
\section{Nondegenerate Markov chains}
\label{secnondeg}

The nondegeneracy Assumption~\ref{asptnondeg} will play an essential role
in our theory. Before we can turn to the proofs of our main results, we
must therefore begin by establishing some general consequences of the
nondegeneracy assumption that will be needed throughout the paper.

%s3.1 #&#
\subsection{Product structure of the invariant measure}

Assumption~\ref{asptnondeg} states that the transition kernel $P$ of
the Markov chain $(X_n,Y_n)_{n\in\mathbb{Z}}$ is equivalent to a product\vadjust{\goodbreak}
of transition kernels of two independent Markov chains. Our first
question is, therefore, whether this forces the invariant measure $\pi$ to
possess a similar product structure; that is, if a stationary Markov
chain is nondegenerate, then is its invariant measure necessarily
equivalent to the product of its marginals? In general, of course, the
answer is negative (e.g., consider the case where $P$ is the
identity and $\pi$ is any probability measure that is not equivalent to a
product measure). However, we will presently show that if, in addition to
nondegeneracy, we assume that the marginal process
$(X_n)_{n\in\mathbb{Z}}$ is ergodic in a suitable sense, then $\pi$ is
forced to possess the desired product structure.

We need two lemmas. The first states that the nondegeneracy of the
transition kernel $P$ implies that the iterates $P^n$ are also
nondegenerate; in fact, we show that
$\mathbf{P}((X_n,Y_n)\in\cdv X_0,Y_0)\sim
\mathbf{P}(X_n\in\cdv X_0)\otimes\mathbf{P}(Y_n\in\cdv Y_0)$.
%
%le3.1 #&#
\begin{lem}
\label{lemnondegexact}
Suppose that Assumption~\ref{asptnondeg} is in force. Choose fixed
versions $\pi^Y(w,dz)$, $\pi^X(z,dw)$ of the regular conditional
probabilities $\mathbf{P}(X_0\in\cdv Y_0)$,
$\mathbf{P}(Y_0\in\cdv X_0)$, respectively, and define the
probability kernels\label{page16}
%(Lemma~\ref{lemnondegexact})}%
%(Lemma~\ref{lemnondegexact})}%
%(Lemma~\ref{lemnondegexact})}%
%(Lemma~\ref{lemnondegexact})}%
%
\begin{eqnarray*}
P_n^X(z,A) &=& \int\mathbf{1}_A(z') P^n(z,w,dz',dw')
\pi^X(z,dw),\\
P_n^Y(w,B) &=& \int\mathbf{1}_B(w') P^n(z,w,dz',dw')
\pi^Y(w,dz).
\end{eqnarray*}
Then we have for all $n\in\mathbb{N}$
\[
P^n(z,w,dz',dw') = G_n(z,w,z',w')
P_n^X(z,dz')P_n^Y(w,dw'),
\]
where $G_n\dvtx E\times F\times E\times F\to\ ]0,\infty[$ are
strictly positive measurable functions.
\end{lem}
\begin{pf}
From the Assumption~\ref{asptnondeg}, it follows directly that
\[
P^n(z,w,dz',dw') = g_n(z,w,z',w')P_0^n(z,dz')Q^n(w,dw')
\]
for some strictly positive measurable function $g_n\dvtx E\times F\times
E\times F\to\ ]0,\infty[$. But then the result follows
directly from the definition of $P_n^X$, $P_n^Y$ with
\begin{eqnarray*}
&&G_n(z,w,z',w') \\
&&\qquad=g_n(z,w,z',w')\\
&&\qquad\quad{}\times\biggl(
\int g_n(z,\tilde w,z',\tilde w')Q^n(\tilde w,d\tilde w')
\pi^X(z,d\tilde w)\\
&&\hspace*{50pt}{}\times\int g_n(\tilde z,w,\tilde z',w')P_0^n(\tilde z,d\tilde z')
\pi^Y(w,d\tilde z)\biggr)^{-1}.
\end{eqnarray*}
The proof is complete.
\end{pf}

The second lemma states that if the unobserved process
$(X_n)_{n\in\mathbb{Z}}$ is ergodic in a suitable sense,
and if the nondegeneracy assumption holds, then every $P$-invariant
function is independent of its unobserved component.
%
%le3.2 #&#
\begin{lem}
\label{lembirkhoffprep}
Suppose that Assumption~\ref{asptnondeg} is in force, and that
\[
\int\|P_n^X(z, \cdot)-\pi( \cdot\times F)\|_{\mathrm{TV}}
\pi(dz\times F) \xrightarrow{n\to\infty}0.
\]
Then for any bounded measurable function $f\dvtx E\times F\to\mathbb{R}$ that
is $P$-invariant (i.e., $f=Pf)$, there exists a bounded measurable
function $g\dvtx F\to\mathbb{R}$ such that $f(z,w)=g(w)$ for $\pi$-a.e.
$(z,w)\in E\times F$.
\end{lem}
\begin{pf}
As $f$ is $P$-invariant, the process $(f(X_n,Y_n))_{n\ge0}$ is a
martingale under~$\mathbf{P}$. By stationarity and the martingale
convergence theorem,
\[
\mathbf{E}\bigl(|f(X_n,Y_n)-f(X_0,Y_0)|\bigr) =
\mathbf{E}\bigl(|f(X_{n+k},Y_{n+k})-f(X_k,Y_k)|\bigr)
\xrightarrow{k\to\infty}0.
\]
In particular, we have
\[
\int\mathbf{P}^{z,w}\bigl(f(X_0,Y_0)=f(X_n,Y_n)\mbox{ for all }n\ge0\bigr)
\pi(dz,dw) = 1.
\]
Therefore, we may choose a set $H_1\in\mathcal{B}(E\times F)$ with
$\pi(H_1)=1$ such that
\[
\mathbf{P}^{z,w}\bigl(f(z,w)=f(X_n,Y_n)\mbox{ for all }n\ge0\bigr) =1
\qquad\mbox{for all }(z,w)\in H_1.
\]
Next, let $\rho\dvtx F\times\mathcal{B}(E)\to[0,1]$ be a version of the regular
conditional probability $\mathbf{P}(X_0\in\cdv Y_0)$. Then
by our assumption and the triangle inequality,
\begin{eqnarray*}
&&\int\|P_n^X(z, \cdot)-P_n^X(z', \cdot)\|_{\mathrm{TV}}
\rho(w,dz)\rho(w,dz')\pi(E\times dw) \\
&&\qquad\le2\int\|P_n^X(z, \cdot)-\pi( \cdot\times F)\|_{\mathrm{TV}}
\pi(dz\times F) \xrightarrow{n\to\infty}0.
\end{eqnarray*}
Therefore, using Fatou's lemma, we can choose a set
$H_2\in\mathcal{B}(E\times E\times F)$ of
$(\rho\otimes\rho)\pi(E\times\cdot)$-full measure such that
\[
{\liminf_{n\to\infty}}\|P_n^X(z, \cdot)-P_n^X(z', \cdot)\|_{\mathrm{TV}}
=0 \qquad\mbox{for all }(z,z',w)\in H_2.
\]
Now define the set $H\in\mathcal{B}(E\times E\times F)$ as follows:
\[
H = \{(z,z',w)\in E\times E\times F\dvtx
(z,w),(z',w)\in H_1\}\cap H_2.
\]
Then it is easily seen that the set $H$ has
$(\rho\otimes\rho)\pi(E\times\cdot)$-full measure.

We now claim that $f(z,w)=f(z',w)$ for every $(z,z',w)\in H$.
To see this, let us fix some point $(z,z',w)\in H$, and choose
$n\ge0$ such that
\[
\|P_n^X(z, \cdot)-P_n^X(z', \cdot)\|_{\mathrm{TV}}<1.
\]
Thus $P_n^X(z, \cdot)$ and $P_n^X(z', \cdot)$ are not mutually
singular. By Lemma~\ref{lemnondegexact}
\[
P^n(z,w, \cdot)\sim P_n^X(z, \cdot)\otimes
P_n^Y(w, \cdot),\qquad
P^n(z',w, \cdot)\sim P_n^X(z', \cdot)\otimes
P_n^Y(w, \cdot).
\]
Therefore, $P^n(z,w, \cdot)$ and $P^n(z',w, \cdot)$ are not
mutually singular. But note that, by the definition of $H$,
$P^n(z,w, \cdot)$ is supported on the set
\[
\Xi_1=\{(\tilde z,\tilde w)\in E\times F\dvtx f(z,w)=
f(\tilde z,\tilde w)\},\vspace*{-1pt}
\]
while $P^n(z',w, \cdot)$ is supported on the set
\[
\Xi_2=\{(\tilde z,\tilde w)\in E\times F\dvtx f(z',w)=
f(\tilde z,\tilde w)\}.\vspace*{-1pt}
\]
Thus the fact that $P^n(z,w, \cdot)$ and $P^n(z',w, \cdot)$ are not
mutually singular implies that $\Xi_1\cap\Xi_2\ne\varnothing$, which
establishes the claim.

To complete the proof, define $g(w) = \int f(z,w) \rho(w,dz)$.
Then
\begin{eqnarray*}
&&\int|f(z,w)-g(w)| \pi(dz,dw) \\[-2pt]
&&\qquad\le\int|f(z,w)-f(z',w)| \rho(w,dz)
\rho(w,dz')\pi(E\times dw) = 0.\vspace*{-1pt}
\end{eqnarray*}
Thus $f(z,w)=g(w)$ for $\pi$-a.e. $(z,w)\in E\times F$ as
desired.\vspace*{-3pt}
\end{pf}

We can now prove the main result of this subsection: if the nondegeneracy
assumption holds, and if, in addition, the unobserved component
$(X_n)_{n\in\mathbb{Z}}$ is ergodic, then the invariant measure $\pi$ is
necessarily equivalent to the product of its marginals. Note that the
ergodicity assumption in this result automatically holds when Assumption
\ref{asptergodic} is in force.\vspace*{-3pt}
%
%pr3.3 #&#
\begin{prop}
\label{propproductpi}
Suppose that Assumption~\ref{asptnondeg} is in force, and that
\[
\int\|P_n^X(z, \cdot)-\pi( \cdot\times F)\|_{\mathrm{TV}}
\pi(dz\times F) \xrightarrow{n\to\infty}0.\vspace*{-1pt}
\]
Then there exists a strictly positive
measurable function $h\dvtx E\times F\to\ ]0,\infty[$ such that
$\pi(dz,dw) = h(z,w)\pi(dz\times F)\pi(E\times dw)$.\vspace*{-3pt}
\end{prop}
\begin{pf}
We begin by noting that
\[
\pi(A\times F) = \int\pi(dz\times F) P_n^X(z,A),\qquad
\pi(E\times B) = \int\pi(E\times dw) P_n^Y(w,B)\vspace*{-1pt}
\]
by the invariance of $\pi$. Now let $C\in\mathcal{B}(E\times F)$ be a set
such that $\pi(C)=0$. As $\pi P^n=\pi$, it follows from Lemma
\ref{lemnondegexact} that
\[
\int\mathbf{1}_C(z',w')P_n^X(z,dz')P_n^Y(w,dw')\pi(dz,dw)=0\vspace*{-1pt}
\]
for all $n\in\mathbb{N}$. But note that
\begin{eqnarray*}
&&\int\mathbf{1}_C(z,w)\pi(dz\times F)\pi(E\times dw) \\[-2pt]
&&\qquad =
\int\mathbf{1}_C(z',w')\pi(dz'\times F)
P_n^Y(w,dw')\pi(dz,dw) \\[-2pt]
&&\qquad \le
\int\|P_n^X(z, \cdot)-\pi( \cdot\times F)\|_{\mathrm{TV}}
\pi(dz\times F) \\[-2pt]
&&\qquad\quad{} +
\int\mathbf{1}_C(z',w')P_n^X(z,dz')P_n^Y(w,dw')\pi(dz,dw).\vspace*{-1pt}
\end{eqnarray*}
Letting $n\to\infty$ and using the ergodicity assumption gives
\[
\int\mathbf{1}_C(z,w)\pi(dz\times F)\pi(E\times dw) = 0.
\]
As this holds for any set $C$ such that $\pi(C)=0$, we have evidently
shown that $\pi(dz\times F)\pi(E\times dw)\ll\pi(dz,dw)$.
Conversely, choose a set $C$ such that
\[
\int\mathbf{1}_C(z,w)\pi(dz\times F)\pi(E\times dw) = 0.
\]
Then, by Lemma~\ref{lemnondegexact}, we have
\[
\int P^n(z,w,C)\pi(dz\times F)\pi(E\times dw) = 0
\]
for all $n\in\mathbb{N}$. By the Birkhoff ergodic theorem,
\[
\frac{1}{N}\sum_{n=1}^NP^n(z,w,C)
\xrightarrow{N\to\infty}f(z,w)\qquad
\mbox{for }\pi\mbox{-a.e. }(z,w)\in E\times F,
\]
where $f$ is a $P$-invariant function with $\pi(f)=\pi(C)$.
Moreover, by Lemma~\ref{lembirkhoffprep} we have
$f(z,w)=g(w)$ for $\pi$-a.e. $(z,w)\in E\times F$ for some function $g$.
But as we have already shown that $\pi(dz\times F)\pi(E\times
dw)\ll\pi(dz,dw)$, these statements hold
$\pi(dz\times F)\pi(E\times dw)$-a.e. also. Therefore,
\begin{eqnarray*}
&&0 = \frac{1}{N}\sum_{n=1}^N\int P^n(z,w,C)\pi(dz\times F)
\pi(E\times dw)\\
&&\quad\xrightarrow{N\to\infty} \quad
\int g(w) \pi(E\times dw) =
\int f(z,w) \pi(dz,dw) = \pi(C).
\end{eqnarray*}
As this holds for any $C$ such that $\int\mathbf{1}_C(z,w)\pi(dz\times
F)\pi(E\times dw) = 0$, we evidently have
$\pi(dz,dw)\ll\pi(dz\times F)\pi(E\times dw)$, and the proof is complete.
\end{pf}

%s3.2 #&#
\subsection{Reversed nondegeneracy}

One important consequence of Proposition~\ref{propproductpi} is that, if
the unobserved process $(X_n)_{n\in\mathbb{Z}}$ is ergodic and the
transition kernel $P$ is nondegenerate, then the nondegeneracy Assumption
\ref{asptnondeg} holds also in reverse time (i.e., the backward
transition kernel $P'$ must be nondegenerate also). In particular, this
implies that the Assumptions~\ref{asptergodic}--\ref{asptnondeg} are
invariant under time reversal.
%
%le3.4 #&#
\begin{lem}
\label{lemnondegrv}
Suppose that Assumption~\ref{asptnondeg} is in force, and that
\[
\int\|P_n^X(z, \cdot)-\pi( \cdot\times F)\|_{\mathrm{TV}}
\pi(dz\times F) \xrightarrow{n\to\infty}0.
\]
Then $P'$ is also nondegenerate; that is, there exist transition
probability kernels $P'_0\dvtx E\times\mathcal{B}(E)\to[0,1]$ and
$Q'\dvtx F\times\mathcal{B}(F)\to[0,1]$ such that
\[
P'(z,w,dz',dw')=g'(z,w,z',w')P'_0(z,dz')Q'(w,dw')
\]
for some strictly positive measurable function $g'\dvtx E\times F\times
E\times F\to\ ]0,\infty[$.
\end{lem}
\begin{pf}
Note that by Proposition~\ref{propproductpi}
and Assumption~\ref{asptnondeg}
\[
\mathbf{E}\bigl((X_0,Y_0,X_1,Y_1)\in B\bigr) =
\int_B h(x_0,y_0)
g(x_0,y_0,x_1,y_1)
\rho(dx_0,dx_1)\kappa(dy_0,dy_1),
\]
where $\rho(dx,dx')=\pi(dx\times F)P_0(x,dx')$,
$\kappa(dy,dy')=\pi(E\times dy)Q(y,dy')$, and where $g,h$ are
strictly positive measurable functions. Let us now fix any versions
$r(x_1,dx_0)$ and $k(y_1,dy_0)$ of the regular conditional probabilities
$\rho(X_0\in\cdv X_1)$ and $\kappa(Y_0\in\cdv Y_1)$,
respectively. Then by the Bayes formula,
\[
\mathbf{P}\bigl((X_0,Y_0)\in A|X_1,Y_1\bigr) =
\frac{\int_A h(z,w)
g(z,w,X_1,Y_1)r(X_1,dz)k(Y_1,dw)
}{\int h(z,w) g(z,w,X_1,Y_1)r(X_1,dz)k(Y_1,dw)}.
\]
As $P'$ is a version of
$\mathbf{P}((X_0,Y_0)\in\cdv X_1,Y_1)$, the result follows.
\end{pf}

%s3.3 #&#
\subsection{Equivalence of the observations}

We now turn to a different consequence of the nondegeneracy assumption.
It is\vspace*{1pt} easily seen that when Assumption~\ref{asptnondeg} holds,
the laws of $(Y_0,\ldots,Y_n)$ under $\mathbf{P}^{z,w}$ and
$\mathbf{P}^{z',w}$ are equivalent for any $z,z'\in E$, $w\in F$,
$n<\infty$. That is, the laws of the observed process under different
initializations of the unobserved process are equivalent on any finite
time horizon. To prove our main results, however, we will require such
an equivalence to hold on the \textit{infinite} time horizon. The following
result is therefore of central importance.
%
%pr3.5 #&#
\begin{prop}
\label{propcoupling}
Suppose that Assumption~\ref{asptnondeg} holds.
Let $\xi,\xi'$ be probability measures on $(E,\mathcal{B}(E))$, let
$\eta$ be a probability measure on $(F,\mathcal{B}(F))$ and let
$v\dvtx E\times F\to\ ]0,\infty[$ and $v'\dvtx E\times
F\to\ ]0,\infty[$ be strictly positive measurable functions.
Define the probability measures on $(E\times F,\mathcal{B}(E\times F))$
\[
\nu(dx,dy) = v(x,y)\xi(dx)\eta(dy),\qquad
\nu'(dx,dy) = v'(x,y)\xi'(dx)\eta(dy).
\]
If
${\liminf_{n\to\infty}}\|\mathbf{P}^\nu(X_n\in\cdot)-\mathbf{P}^{\nu
'}(X_n\in\cdot)
\|_{\mathrm{TV}}=0$, then
$\mathbf{P}^\nu|_{\mathcal{F}^Y_+}\sim
\mathbf{P}^{\nu'}|_{\mathcal{F}^Y_+}$.
\end{prop}
\begin{pf}
Choose any $A\in\mathcal{F}^Y_+$ such that
$\mathbf{P}^{\nu'}(A)=0$. It suffices to prove that
$\mathbf{P}^\nu(A)=0$. Indeed, this shows that
$\mathbf{P}^\nu|_{\mathcal{F}^Y_+}\ll
\mathbf{P}^{\nu'}|_{\mathcal{F}^Y_+}$, while the reverse
statement follows as the assumptions are symmetric in $\nu$ and
$\nu'$.

Fix for the time being $n\in\mathbb{N}$. Note that by construction
\[
I_A(x,y) = I_A(y(0),\ldots,y(n),(y(k))_{k>n}).\vadjust{\goodbreak}
\]
Let us define the measurable function
\[
a(y_0,\ldots,y_n,x_n) = \mathbf{E}^{x_n,y_n}(
I_A(y_0,\ldots,y_{n},(Y_k)_{k\ge1})).
\]
Then, by the Markov property,
\[
a(Y_0,\ldots,Y_n,X_n) =
\mathbf{P}^\rho\bigl(A|\mathcal{F}^X_{[0,n]}\vee\mathcal{F}^Y_{[0,n]}\bigr),\qquad
\mathbf{P}^\rho\mbox{-a.s.}
\]
for any initial probability measure $\rho$. In particular,
\[
a(Y_0,\ldots,Y_n,X_n) = 0 ,\qquad\mathbf{P}^{\nu'}\mbox{-a.s.}
\]
Let $\mathbf{Q}^\eta$ be the law of the Markov chain $(Y_k)_{k\ge
0}$ with initial measure $\eta$ and transition kernel $Q$, and
$\mathbf{P}_0^{\xi'}$ be the law of the Markov chain $(X_k)_{k\ge
0}$ with initial measure $\xi'$ and transition kernel $P_0$.
By our assumptions,
\[
\mathbf{P}^{\nu'}(A)=
(\mathbf{P}_0^{\xi'}\otimes\mathbf{Q}^\eta)
\Biggl[I_A
v'(X_0,Y_0)\prod_{i=0}^{n-1}g(X_i,Y_i,X_{i+1},Y_{i+1})
\Biggr]
\]
for every $A\in\mathcal{F}^X_{[0,n]}\vee\mathcal{F}^Y_{[0,n]}$.
In particular, the law of
$(Y_0,\ldots,Y_n,X_n)$ under~$\mathbf{P}^{\nu'}$ and the law of
$\mathbf{Q}^\eta|_{\mathcal{F}^Y_{[0,n]}}\otimes\xi'P_0^n$ are
equivalent. Therefore,
\[
a(Y_0,\ldots,Y_n,X_n)=0,\qquad
\bigl(\mathbf{Q}^\eta|_{\mathcal{F}^Y_{[0,n]}}\otimes
\xi'P_0^n\bigr)\mbox{-a.s.}
\]
Choose $S_n\in\mathcal{B}(E)$ such that $(\xi
P_0^n)( \cdot\cap S_n)\ll\xi' P_0^n$ and $(\xi'
P_0^n)(S_n^c)=0$ (so~$S_n$ defines the Lebesgue decomposition of
$\xi P_0^n$ with respect to $\xi' P_0^n$). Then
\[
I_{S_n}(X_n) a(Y_0,\ldots,Y_n,X_n)=0,\qquad
\bigl(\mathbf{Q}^\eta|_{\mathcal{F}^Y_{[0,n]}}\otimes
\xi P_0^n\bigr)\mbox{-a.s.}
\]
Therefore,
\[
a(Y_0,\ldots,Y_n,X_n) \le I_{S_n^c}(X_n),\qquad
\bigl(\mathbf{Q}^\eta|_{\mathcal{F}^Y_{[0,n]}}\otimes
\xi P_0^n\bigr)\mbox{-a.s.}
\]
But, as above, we find that the law of $(Y_0,\ldots,Y_n,X_n)$ under
$\mathbf{P}^{\nu}$ is equivalent to
$\mathbf{Q}^\eta|_{\mathcal{F}^Y_{[0,n]}}\otimes\xi P_0^n$. Therefore,
we obtain immediately
\[
a(Y_0,\ldots,Y_n,X_n) \le I_{S_n^c}(X_n),\qquad
\mathbf{P}^{\nu}\mbox{-a.s.}
\]
Taking the expectation, we find that $\mathbf{P}^\nu(A) \le
\mathbf{P}^\nu(X_n\in S_n^c)$.

At this point, we note that $n\in\mathbb{N}$ in the above
construction was arbitrary.
Moreover, we have already shown that for any $n\in\mathbb{N}$, the law of
$X_n$ under~$\mathbf{P}^{\nu'}$ is equivalent to $\xi' P_0^n$.
Therefore $\mathbf{P}^{\nu'}(X_n\in S_n^c)=0$, and we find that
\[
\mathbf{P}^\nu(A)\le\liminf_{n\to\infty}
\mathbf{P}^\nu(X_n\in S_n^c) \le
{\liminf_{n\to\infty}}
\|\mathbf{P}^\nu(X_n\in\cdot)-
\mathbf{P}^{\nu'}(X_n\in\cdot)\|_{\mathrm{TV}}
=0.
\]
Thus the proof is complete.
\end{pf}

A useful corollary is the following result.
%
%co3.6 #&#
\begin{cor}
\label{coryinfhabscont} Suppose that Assumptions
\ref{asptergodic} and~\ref{asptnondeg} hold, and let $\mu$ be a
probability measure on $E\times F$ such that
$\mu(E\times\cdot)\ll\pi(E\times\cdot)$ and
\[
\mathbf{E}^\mu\bigl(
\bigl\|\mathbf{P}^\mu(X_n\in\cdv Y_0)-\mathbf{P}(X_n\in\cdot)
\bigr\|_{\mathrm{TV}}\bigr)
\xrightarrow{n\to\infty}0.
\]
Then $\mathbf{P}^\mu|_{\mathcal{F}^Y_+}\ll\mathbf{P}|_{\mathcal{F}^Y_+}$.
If $\mu(E\times\cdot)\sim\pi(E\times\cdot)$, then
$\mathbf{P}^\mu|_{\mathcal{F}^Y_+}\sim\mathbf{P}|_{\mathcal{F}^Y_+}$.
\end{cor}
\begin{pf}
We begin by noting that
\[
\mathbf{E}\bigl(
\bigl\|\mathbf{P}(X_n\in\cdv Y_0)-\mathbf{P}(X_n\in\cdot)
\bigr\|_{\mathrm{TV}}\bigr) \le
\mathbf{E}\bigl(
\|\mathbf{P}^{X_0,Y_0}(X_n\in\cdot)-\mathbf{P}(X_n\in\cdot)
\|_{\mathrm{TV}}\bigr).
\]
Therefore, by Assumption~\ref{asptergodic},
\[
\bigl\|\mathbf{P}(X_n\in\cdv Y_0)-\mathbf{P}(X_n\in\cdot)
\bigr\|_{\mathrm{TV}}\xrightarrow{n\to\infty}0
\qquad\mbox{in }\mathbf{P}\mbox{-probability}.
\]
As $\mathbf{P}^\mu(Y_0\in\cdot)\ll\mathbf{P}(Y_0\in\cdot)$,
this convergence is also in $\mathbf{P}^\mu$-probability. Therefore,
using dominated convergence and the triangle inequality,
\[
\mathbf{E}^\mu\bigl(
\bigl\|\mathbf{P}^\mu(X_n\in\cdv Y_0)-\mathbf{P}(X_n\in\cdv Y_0)
\bigr\|_{\mathrm{TV}}\bigr)
\xrightarrow{n\to\infty}0.
\]
By Fatou's lemma, we obtain
\[
{\liminf_{n\to\infty}}
\bigl\|\mathbf{P}^\mu(X_n\in\cdv Y_0)-\mathbf{P}(X_n\in\cdv Y_0)
\bigr\|_{\mathrm{TV}}=0,\qquad
\mathbf{P}^\mu\mbox{-a.s.}
\]
Let $\nu\dvtx F\times\mathcal{B}(E)\to[0,1]$,
$\nu'\dvtx F\times\mathcal{B}(E)\to[0,1]$ be versions of the regular
conditional probabilities $\mathbf{P}^\mu(X_0\in\cdv Y_0)$,
$\mathbf{P}(X_0\in\cdv Y_0)$, respectively. Then
\[
\liminf_{n\to\infty}
\bigl\|
\mathbf{P}^{\nu(w, \cdot)\otimes\delta_w}(X_n\in\cdot) -
\mathbf{P}^{\nu'(w, \cdot)\otimes\delta_w}(X_n\in\cdot)
\bigr\|_{\mathrm{TV}}=0,\qquad \mu(E\times\cdot)
\mbox{-a.e. }w.
\]
By Proposition~\ref{propcoupling}, it follows that
\[
\mathbf{P}^{\nu(w, \cdot)\otimes\delta_w}|_{\mathcal{F}^Y_+}\sim
\mathbf{P}^{\nu'(w, \cdot)\otimes\delta_w}|_{\mathcal{F}^Y_+},\qquad
\mu(E\times\cdot)\mbox{-a.e. }w.
\]
By the Lebesgue decomposition for kernels (\cite{DM}, Section V.58),
there is a~measurable version of the Radon--Nikodym derivative.
It follows that
\[
\mathbf{P}^\mu|_{\mathcal{F}^Y_+}
=
\mathbf{P}^{\nu\otimes\mu(E\times\cdot)}|_{\mathcal{F}^Y_+}
\sim
\mathbf{P}^{\nu'\otimes\mu(E\times\cdot)}|_{\mathcal{F}^Y_+}
\ll
\mathbf{P}^{\nu'\otimes\pi(E\times\cdot)}|_{\mathcal{F}^Y_+}
=
\mathbf{P}|_{\mathcal{F}^Y_+},
\]
where we have used that $\mu(E\times\cdot)\ll\pi(E\times\cdot)$. If
$\mu(E\times\cdot)\sim\pi(E\times\cdot)$, then clearly $\ll$
can be replaced by $\sim$ in the previous equation.
\end{pf}

%s4 #&#
\section{\texorpdfstring{Proof of Theorem \protect\ref{thmexchg}}{Proof of Theorem 2.9}}
\label{secexchg}

The goal of this section is to prove Theorem~\ref{thmexchg}. To this end,
we begin by recalling the basic result from~\cite{vH10} on the ergodicity
of Markov chains in random environments. This result will be used to
establish that the unobservable process $(X_n)_{n\ge0}$ has
trivial tail $\sigma$-field under the conditional measure
$\mathbf{P}( \cdv \mathcal{F}^Y)$. Finally, we show that
$\mathbf{P}(X_0\in\cdv \mathcal{F}^Y)\sim
\mathbf{P}(X_0\in\cdv \mathcal{F}^Y_+)$, which allows us to complete the
proof by applying a result of von Weizs\"acker~\cite{vW83}.

%s4.1 #&#
\subsection{Markov chains in random environments}
\label{secmcre}

We begin by recalling the relevant notions from~\cite{vH10}, Section 2.
A Markov chain in a random environment is defined by the following three
ingredients:
\begin{longlist}[(2)]
\item[(1)] A probability kernel $P^X\dvtx E\times\Omega^Y\times
\mathcal{B}(E)\to[0,1]$.
\item[(2)] A probability kernel $\varpi\dvtx \Omega^Y\times\mathcal{B}(E)\to[0,1]$
such that
\[
\int P^X(z,y,A) \varpi(y,dz) = \varpi(\Theta y,A)\qquad
\mbox{for all } y\in\Omega^Y, A\in\mathcal{B}(E).
\]
\item[(3)] A stationary probability measure $\mathbf{P}^Y$ on
$(\Omega^Y,\mathcal{F}^Y)$.
\end{longlist}
The process $X_n$ is called a Markov chain in a random environment when
\begin{eqnarray*}
P^X(X_n,Y\circ\Theta^n,A) &=& \mathbf{P}(X_{n+1}\in A|
\mathcal{F}_n^X\vee\mathcal{F}^Y),\qquad
\mathbf{P}\mbox{-a.s.},\\
\varpi(Y\circ\Theta^n,A) &=& \mathbf{P}(X_n\in A|\mathcal{F}^Y),\qquad
\mathbf{P}\mbox{-a.s.}
\end{eqnarray*}
for every $A\in\mathcal{B}(E)$ and $n\in\mathbb{Z}$, and
$\mathbf{P}^Y = \mathbf{P}|_{\mathcal{F}^Y}$.\label{page24}
One should think of a~Markov chain in a random environment $X_n$ as a
process that is Markov \textit{conditionally} on the environment
$Y$. The conditional chain is time-inhomogeneous but must satisfy
certain stationarity properties: the environment is stationary and the
(time-dependent) conditional transition probabilities
\mbox{$P^X( \cdot,Y\circ\Theta^n, \cdot)$} and quasi-invariant measure
$\varpi(Y\circ\Theta^n, \cdot)$ are themselves stationary processes
with respect to the environment. The stationarity properties ensure
that Markov chains in random environments behave ``almost'' like
time-homogeneous Markov chains; cf. Theorem~\ref{thmmcre} below.

Let us introduce a probability kernel
$\mathbf{P}_{\cdot}\dvtx E\times\Omega^Y\times\mathcal{F}^X_+\to[0,1]$
so that
\begin{eqnarray*}
\mathbf{P}_{z,y}(A) &=&
\int I_A(x) P^X\bigl(x(n-1),\Theta^{n-1}y,dx(n)\bigr)\times \cdots\\
&&\hspace*{9.5pt}{}\times
P^X(x(1),\Theta y,dx(2)) P^X(x(0),y,dx(1)) \delta_z(dx(0))
\end{eqnarray*}
for $A\in\mathcal{F}_{[0,n]}^X$. It is easily seen that
$\mathbf{P}_{z,y}$ is a version of the regular conditional probability
$\mathbf{P}((X_k)_{k\ge
0}\in\cdv \mathcal{F}_0^X\vee\mathcal{F}^Y)$.\label{pageCondlaw}
%given $X_0=z$, $Y=y$}%
We can now state the following ergodic theorem for Markov chains in random
environments (\cite{vH10}, Theorem 2.3).
%
%th4.1 #&#
\begin{theorem}
\label{thmmcre}
The following are equivalent.
\begin{longlist}[(1)]
\item[(1)] $\|\mathbf{P}_{z,y}(X_n\in\cdot)-
\mathbf{P}_{z',y}(X_n\in\cdot)\|_{\mathrm{TV}}\xrightarrow{n\to\infty}0$
for $(\varpi\otimes\varpi)\mathbf{P}^Y$-a.e. $(z,z',y)$.
\item[(2)] The tail $\sigma$-field
$\mathcal{T}^X=\bigcap_{n\ge0}\mathcal{F}^X_{[n,\infty[}$
is a.s. trivial in the following sense:
\[
\mathbf{P}_{z,y}(A)=\mathbf{P}_{z,y}(A)^2=
\mathbf{P}_{z',y}(A)\qquad
\mbox{for all }
A\in\mathcal{T}^X\mbox{ and }(z,z',y)\in H,
\]
where $H$ is a fixed set (independent of $A$) of
$(\varpi\otimes\varpi)\mathbf{P}^Y$-full measure.
\item[(3)] For $(\varpi\otimes\varpi)\mathbf{P}^Y$-a.e. $(z,z',y)$, there
is an $n\in\mathbb{N}$ such that the measures
$\mathbf{P}_{z,y}(X_n\in\cdot)$ and
$\mathbf{P}_{z',y}(X_n\in\cdot)$
are not mutually singular.
\end{longlist}
\end{theorem}

%s4.2 #&#
\subsection{Weak ergodicity of the conditional process}

Our first order of business is to establish that, under the model defined
in this paper, $X_n$ is indeed a Markov chain in a random environment in
the sense of Section~\ref{secmcre}, where the observations $Y$ play the
role of the environment; that is, we must show that the unobserved process
$X_n$ is still a Markov chain \textit{conditionally} on the observations $Y$
satisfying the requisite stationarity properties. This is the statement of
the following lemma, whose proof is omitted as it is identical to that in
\cite{vH10}. As everything that follows is based on this elementary fact,
however, let us briefly sketch why the result is true for the convenience
of the reader. It is easily seen that
\[
\mathbf{P}(X_{n+1}\in\cdv
\mathcal{F}_n^X\vee\mathcal{F}^Y)\circ\Theta^{-n} =
\mathbf{P}(X_{1}\in\cdv
\mathcal{F}_0^X\vee\mathcal{F}^Y) =
\mathbf{P}(X_{1}\in\cdv
\sigma\{X_0\}\vee\mathcal{F}^Y).
\]
The first equality follows from stationarity of $\mathbf{P}$, and the
second equality follows as $\mathcal{F}_{[1,\infty[}$ is conditionally
independent of $\mathcal{F}_{-1}$ given $\sigma\{X_0,Y_0\}$
by the Markov property of $(X_n,Y_n)_{n\in\mathbb{Z}}$. We can therefore
choose $P^X$ to be a regular version of
$\mathbf{P}(X_{1}\in\cdv \sigma\{X_0\}\vee\mathcal{F}^Y)$.
Similarly, we can choose $\varpi$ to be a regular version of
$\mathbf{P}(X_0\in\cdv \mathcal{F}^Y)$, and $\mathbf{P}^Y$ to be the
law of $Y$. It is now an elementary exercise to check that these kernels
do indeed characterize the process $X_n$ as a Markov chain in a random
environment in the sense of Section~\ref{secmcre}.
%
%le4.2 #&#
\begin{lem}
\label{lemhmmismcre}
There exist probability kernels $P^X\dvtx E\times\Omega^Y\times
\mathcal{B}(E)\to[0,1]$ and $\varpi\dvtx\Omega^Y\times\mathcal{B}(E)\to[0,1]$,
and a probability measure $\mathbf{P}^Y$ on $(\Omega^Y,\mathcal{F}^Y)$,
such that the conditions of Section~\ref{secmcre} are satisfied.
\end{lem}
\begin{pf}
The proof is identical to that of~\cite{vH10}, Lemma 3.3.
\end{pf}

The main goal of this subsection is to prove the following theorem.
%
%th4.3 #&#
\begin{theorem}
\label{thmerg} Suppose that both Assumptions~\ref{asptergodic} and
\ref{asptnondeg} are in force. Then any (hence, all) of the
conditions of Theorem~\ref{thmmcre} hold true.
\end{theorem}

The strategy of the proof of Theorem~\ref{thmerg} is to show that
condition (3) of Theorem~\ref{thmmcre} follows from Assumptions
\ref{asptergodic} and~\ref{asptnondeg}. To this end, we begin by
proving that Theorem~\ref{thmerg} would follow if we can establish
equivalence of the conditional and unconditional transition kernels $P^X$
and $P$.
%
%le4.4 #&#
\begin{lem}
\label{lemwergabsct} Suppose that Assumptions~\ref{asptergodic} and
\ref{asptnondeg} are in force, and that there exists a strictly positive
measurable function $h\dvtx E\times\Omega^Y\times
E\to\ ]0,\infty[$ such that
\[
P^X(z,y,A) =
\int_{E\times F} I_A(\tilde z) h(z,y,\tilde z) P(z,y(0),d\tilde
z,d\tilde{w})\qquad
\mbox{for all }A\in\mathcal{B}(E)
\]
for $\varpi\mathbf{P}^Y$-a.e. $(z,y)$.
Then condition $3$ of Theorem~\ref{thmmcre} holds.
\end{lem}
\begin{pf}
By Assumption~\ref{asptergodic} and the triangle inequality
\[
\int\bigl\|\mathbf{P}^{z,y(0)}(X_n\in\cdot)-
\mathbf{P}^{z',y(0)}(X_n\in\cdot)\bigr\|_{\mathrm{TV}}
\varpi(y,dz)\varpi(y,dz')\mathbf{P}^Y(dy)
\xrightarrow{n\to\infty}0.
\]
By Fatou's lemma, there is a set $H_1$ of
$(\varpi\otimes\varpi)\mathbf{P}^Y$-full measure such that
\[
\liminf_{n\to\infty}\bigl\|\mathbf{P}^{z,y(0)}(X_n\in\cdot)
-\mathbf{P}^{z',y(0)}(X_n\in\cdot)\bigr\|_{\mathrm{TV}}
=0 \qquad\mbox{for all }(z,z',y)\in H_1.
\]
In particular, there is for every $(z,z',y)\in H_1$ an
$n\in\mathbb{N}$ such that $\mathbf{P}^{z,y(0)}(X_n\in\cdot)$ and
$\mathbf{P}^{z',y(0)}(X_n\in \cdot)$ are not mutually singular.

Now let $H_2$ be a set of $\varpi\mathbf{P}^Y$-full measure such that
the absolute continuity condition in the statement of the lemma
holds true for all $(z,y)\in H_2$. By Lem\-ma~\ref{lemmcreinvset},
there is a subset $H_3\subset H_2$ of $\varpi\mathbf{P}^Y$-full measure
such that for every $(z,y)\in H_3$ we have
$\mathbf{P}_{z,y}((X_n,\Theta^ny)\in H_3\mbox{ for all }n\ge0)=1$.
It follows directly that for every $(z,y)\in H_3$, $n\in\mathbb{N}$
and $A\in\mathcal{B}(E)$, we have
\[
\mathbf{P}_{z,y}(X_n\in A)
=\int I_A(x_n)f(x_0,\ldots,x_n,y)\prod_{i=0}^{n-1}
P_0(x_i,dx_{i+1}) \delta_z(dx_0),
\]
where we have defined the strictly positive measurable function
\[
f(x_0,\ldots,x_n,y)=
\prod_{i=0}^{n-1}
h(x_i,\Theta^iy,x_{i+1})\int
g(x_i,y(i),x_{i+1},\tilde{w})
Q(y(i),d\tilde{w}).
\]
On the other hand, we have for every $z,y$
\[
\mathbf{P}^{z,y(0)}(X_n\in A) =
\int I_A(x_n)f'(x_0,\ldots,x_n,y(0))\prod_{i=0}^{n-1}
P_0(x_i,dx_{i+1}) \delta_z(dx_0),
\]
where we have defined the strictly positive measurable function
\[
f'(x_0,\ldots,x_n,y_0)=
\int\prod_{i=0}^{n-1}g(x_i,y_i,x_{i+1},y_{i+1})
Q(y_i,dy_{i+1}).
\]
Therefore $\mathbf{P}_{z,y}(X_n\in\cdot)\sim
\mathbf{P}^{z,y(0)}(X_n\in\cdot)$ for all $(z,y)\in H_3$ and
$n\in\mathbb{N}$.

To complete the proof, define the following set:
\[
H_4=\{(z,z',y)\dvtx(z,z',y)\in H_1, (z,y),(z',y)\in H_3\}.
\]
Then $H_4$ has $(\varpi\otimes\varpi)\mathbf{P}^Y$-full measure, and for
every $(z,z',y)\in H_4$, there is an $n\in\mathbb{N}$ such that
$\mathbf{P}_{z,y}(X_n\in\cdot)$ and
$\mathbf{P}_{z',y}(X_n\in\cdot)$ are not mutually singular.
This establishes condition (3) of Theorem~\ref{thmmcre}.
\end{pf}

We now proceed to prove the following lemma, which verifies the
assumption of Lemma~\ref{lemwergabsct}. This completes the proof of
Theorem~\ref{thmerg}.
%
%le4.5 #&#
\begin{lem}
\label{lemabsctcg}
Suppose that Assumptions~\ref{asptergodic} and
\ref{asptnondeg} hold. Then there exists a strictly positive
measurable function $h\dvtx E\times\Omega^Y\times
E\to\ ]0,\infty[$ such that
\[
P^X(z,y,A) =
\int_{E\times F} I_A(\tilde z) h(z,y,\tilde z) P(z,y(0),d\tilde
z,d\tilde{w})\qquad
\mbox{for all }A\in\mathcal{B}(E),
\]
for $\varpi\mathbf{P}^Y$-a.e. $(z,y)$.
\end{lem}
\begin{pf}
By definition, $P^X$ is a version of the regular conditional probability
$\mathbf{P}(X_1\in\cdv \mathcal{F}_0^X\vee\mathcal{F}^Y)$. But by the
Markov property of $(X_n,Y_n)_{n\in\mathbb{Z}}$, the $\sigma$-fields
$\mathcal{F}_{[1,\infty[}$ and $\mathcal{F}_{-1}$ are conditionally
independent given $\sigma(X_0,Y_0)$. Therefore, $P^X$ is, in fact, a
version of the regular conditional probability
$\mathbf{P}(X_1\in\cdv \sigma(X_0,Y_0)\vee\mathcal{F}^Y_{[1,\infty[})$.
Moreover, clearly the kernel $\tilde P$ defined as
\[
\tilde P(z,w,A) =
\int I_A(\tilde z) P(z,w,d\tilde z,d\tilde w)
\qquad\mbox{for all }A\in\mathcal{B}(E), (z,w)\in E\times F
\]
is a version of the regular conditional probability
$\mathbf{P}(X_1\in\cdv \sigma(X_0,Y_0))$. Finally, we fix throughout
the proof arbitrary versions
$R\dvtx E\times F\times\mathcal{F}^Y_+\rightarrow[0,1]$ and
$R^X\dvtx E\times F\times E\times\mathcal{F}^Y_+\rightarrow[0,1]$
of the regular conditional probabilities
$P((Y_k)_{k\geq1}\in\cdv \sigma(X_0,Y_0))$ and
$P((Y_k)_{k\geq1}\in\cdv \sigma(X_0,Y_0,X_1))$, respectively.
To complete the proof, it suffices to show that
$R^X(z,w,z', \cdot )\sim R(z,w, \cdot)$ for
$(z,w,z')\in H$ with $P((X_0,Y_0,X_1)\in H)=1$.
Indeed, if this is the case, then by the Lebesgue decomposition for
kernels (\cite{DM}, Section V.58), there is a~strictly positive
measurable function $h\dvtx E\times\Omega^Y\times E\to\ ]0,\infty[$
such that
\[
R^X(z,y(0),\tilde z,A) =
\int I_A((y(i))_{i\ge1})
h(z,y,\tilde z) R(z,y(0),d(y(i))_{i\ge1})
\]
for all $A\in\mathcal{F}^Y_{[1,\infty[}$ and $(z,y(0),z')\in H'$ with
$P((X_0,Y_0,X_1)\in H')=1$. It remains to apply Lemma
\ref{lempolishbayes} to the law of the triple
$((X_0,Y_0),X_1,(Y_k)_{k\ge1})$.

It therefore remains to show that $R^X(z,w,z', \cdot )\sim
R(z,w, \cdot)$. To this end, let us introduce convenient versions of
the regular conditional probabilities~$R$ and $R^X$. Note that
we can write for $A\in\mathcal{F}^Y_+$
\[
R(X_0,Y_0,A\circ\Theta) =
\mathbf{E}(\mathbf{P}^{X_1,Y_1}(A)|\sigma(X_0,Y_0)) =
\mathbf{P}^{\nu_{X_0,Y_0}}(A)
\]
by the Markov property of $(X_n,Y_n)_{n\ge0}$,
where we have defined
\[
\nu_{z,w}(d\tilde z,d\tilde w) = P(z,w,d\tilde z,d\tilde w)
= g(z,w,\tilde z,\tilde w)P_0(z,d\tilde z)Q(w,d\tilde w).
\]
On the other hand, using the Bayes formula, we can compute for
$A\in\mathcal{F}^Y_+$
\[
R^X(X_0,Y_0,X_1,A\circ\Theta) =
\mathbf{P}(\mathbf{P}^{X_1,Y_1}(A)|\sigma(X_0,Y_0,X_1)) =
\mathbf{P}^{\nu_{X_0,Y_0,X_1}}(A),
\]
where we have defined
\[
\nu_{z,w,z'}(d\tilde z,d\tilde w) =
\frac{g(z,w,z',\tilde w)}{\int g(z,w,z',w')Q(w,dw')}
\delta_{z'}(d\tilde z)Q(w,d\tilde w).
\]
It therefore\vspace*{1pt} suffices to show that
$\mathbf{P}^{\nu_{z,w,z'}}|_{\mathcal{F}^Y_+}\sim
\mathbf{P}^{\nu_{z,w}}|_{\mathcal{F}^Y_+}$ for $(z,w,z')\in H$ with
$\mathbf{P}((X_0,Y_0,X_1)\in H)=1$. By Proposition~\ref{propcoupling},
it suffices to show that
\[
{\liminf_{n\to\infty}}\|\mathbf{P}^{\nu_{z,w,z'}}(X_n\in\cdot)-
\mathbf{P}^{\nu_{z,w}}(X_n\in\cdot)\|_{\mathrm{TV}}=0
\]
for $(z,w,z')\in H$ with $\mathbf{P}((X_0,Y_0,X_1)\in H)=1$. But
\begin{eqnarray*}
&&\mathbf{E}\bigl(\|\mathbf{P}^{\nu_{X_0,Y_0,X_1}}(X_n\in\cdot)-
\mathbf{P}^{\nu_{X_0,Y_0}}(X_n\in\cdot)\|_{\mathrm{TV}}\bigr) \\
&&\qquad =
\mathbf{E}\bigl(\bigl\|\mathbf{P}(X_{n+1}\in\cdv X_0,Y_0,X_1)-
\mathbf{P}(X_{n+1}\in\cdv X_0,Y_0)\bigr\|_{\mathrm{TV}}\bigr) \\
&&\qquad \le
\mathbf{E}\bigl(\bigl\|\mathbf{P}(X_{n+1}\in\cdv X_0,Y_0,X_1,Y_1)-
\mathbf{P}(X_{n+1}\in\cdot)\bigr\|_{\mathrm{TV}}\bigr) \\
&&\qquad\quad{} +
\mathbf{E}\bigl(\bigl\|\mathbf{P}(X_{n+1}\in\cdv X_0,Y_0)-
\mathbf{P}(X_{n+1}\in\cdot)\bigr\|_{\mathrm{TV}}\bigr) \\
&&\qquad =
\mathbf{E}\bigl(\bigl\|\mathbf{P}(X_{n}\in\cdv X_0,Y_0)-
\mathbf{P}(X_{n}\in\cdot)\bigr\|_{\mathrm{TV}}\bigr) \\
&&\qquad\quad{} +
\mathbf{E}\bigl(\bigl\|\mathbf{P}(X_{n+1}\in\cdv X_0,Y_0)-
\mathbf{P}(X_{n+1}\in\cdot)\bigr\|_{\mathrm{TV}}\bigr),
\end{eqnarray*}
where we have used the triangle inequality and the stationarity of
$\mathbf{P}$. Thus the result follows from Assumption~\ref{asptergodic}
and Fatou's lemma.
\end{pf}

%s4.3 #&#
\subsection{Exchange of intersection and supremum of $\sigma$-fields}

Fix a version $\varpi^+\dvtx\Omega^Y\times\mathcal{B}(E)\to[0,1]$ of the
regular conditional probability
$\mathbf{P}(X_0\in\cdv \mathcal{F}^Y_+)$.\label{page28}
%0}$}%
We begin by establishing the validity of the exchange of intersection and
supremum in Theorem~\ref{thmexchg} assuming that $\varpi^+$ has a
positive density with respect to~$\varpi$.
%
%pr4.6 #&#
\begin{prop}
\label{propvw}
Suppose Assumptions~\ref{asptergodic} and~\ref{asptnondeg} hold,
and that there exists a strictly positive measurable function
$k\dvtx\Omega^Y\times E\to\ ]0,\infty[$ such that
\[
\varpi(y,A) = \int I_A(z)k(y,z)\varpi^+(y,dz)\qquad
\mbox{for all }A\in\mathcal{B}(E)
\]
for $\mathbf{P}^Y$-a.e. $y\in\Omega^Y$.
Then $\bigcap_{n\ge0}\mathcal{F}^Y_+\vee\mathcal{F}^X_{[n,\infty[}
= \mathcal{F}^Y_+$ $\mathbf{P}$-a.s.
\end{prop}
\begin{pf}
By Theorem~\ref{thmerg}, there is a set $H$ of
$(\varpi\otimes\varpi)\mathbf{P}^Y$-full measure with
\[
\mathbf{P}_{z,y}(A)=\mathbf{P}_{z,y}(A)^2=
\mathbf{P}_{z',y}(A)\qquad
\mbox{for all }
A\in\mathcal{T}^X\mbox{ and }(z,z',y)\in H.
\]
As $H$ has $(\varpi\otimes\varpi)\mathbf{P}^Y$-full measure, there clearly
exists a set $H^Y\in\mathcal{B}(\Omega^Y)$ of $\mathbf{P}^Y$-full
measure such that $\int I_H(z,z',y)\varpi(y,dz)\varpi(y,dz')=1$ for all
$y\in H^Y$. Let us now define $\mathbf{P}_y(A)=\int
\mathbf{P}_{z,y}(A)\varpi(y,dz)$. Then
\[
\mathbf{P}_y(A)-\mathbf{P}_y(A)^2 =
\int I_H(z,z',y)\mathbf{P}_{z,y}(A)\bigl(1-\mathbf{P}_{z',y}(A)\bigr)
\varpi(y,dz)\varpi(y,dz') = 0
\]
for every $y\in H^Y$ and $A\in\mathcal{T}^X$. Thus $\mathcal{T}^X$ is
$\mathbf{P}_y$-trivial for all $y\in H^Y$. Therefore, defining
$\mathbf{P}_y^+(A)=\int\mathbf{P}_{z,y}(A)\varpi^+(y,dz)$, our assumption
that $\varpi^+(y, \cdot)\sim\varpi(y, \cdot)$ $\mathbf{P}^Y$-a.e.
$y\in\Omega^Y$ implies that $\mathcal{T}^X$
is $\mathbf{P}_y^+$-trivial for $\mathbf{P}^Y$-a.e.
$y\in\Omega^Y$.\vadjust{\goodbreak}

Now recall that, by definition, $\mathbf{P}_{z,y}$ is a version of the
regular conditional probability $\mathbf{P}((X_k)_{k\ge0}\in\cdv
\sigma(X_0)\vee\mathcal{F}^Y)$. But as $\mathcal{F}^Y_{-1}$ is
conditionally independent of $\mathcal{F}_+$ given $\sigma(X_0,Y_0)$ by
the Markov property, it follows that $\mathbf{P}_{z,y}$ is even a
version of $\mathbf{P}((X_k)_{k\ge0}\in\cdv
\sigma(X_0)\vee\mathcal{F}^Y_+)$. It therefore follows that
$\mathbf{P}_y^+$ is a version of the regular conditional probability
$\mathbf{P}((X_k)_{k\ge0}\in\cdv \mathcal{F}^Y_+)$. We have
therefore shown that $\mathcal{T}^X$ is
$\mathbf{P}( \cdv \mathcal{F}^Y_+)$-trivial $\mathbf{P}$-a.s., which
implies
\[
\bigcap_{n\ge0} \mathcal{F}^Y_+\vee\mathcal{F}^X_{[n,\infty[} =
\mathcal{F}^Y_+,\qquad
\mathbf{P}\mbox{-a.s.}
\]
by Lemma~\ref{lemvWexchg} in Appendix~\ref{secapp}.
\end{pf}

To prove Theorem~\ref{thmexchg}, we must therefore establish that
$\varpi^+$ has a postive density with respect to $\varpi$. It is here
that the time-reversed Assumption~\ref{asptergodicrv} enters the picture;
indeed, the alert reader will not have failed to notice that we have only
used Assumptions~\ref{asptergodic} and~\ref{asptnondeg} up to this point.
%
%le4.7 #&#
\begin{lem}
\label{lemmumuplus}
Suppose that Assumptions~\ref{asptergodic}--\ref{asptnondeg}
are in force. Then there exists a strictly positive measurable function
$k\dvtx\Omega^Y\times E\to\ ]0,\infty[$ such that
\[
\varpi(y,A) = \int I_A(z)k(y,z)\varpi^+(y,dz)\qquad
\mbox{for all }A\in\mathcal{B}(E)
\]
for $\mathbf{P}^Y$-a.e. $y\in\Omega^Y$.
\end{lem}
\begin{pf}
By the Markov property of $(X_n,Y_n)_{n\in\mathbb{Z}}$, we find
that\break
$\mathbf{P}^{z,y_0}((Y_k)_{k<0}\in\cdot)$ and
$\mathbf{P}^{\varpi^+(y, \cdot)\otimes\delta_{y_0}}((Y_k)_{k<0}\in
\cdot)$
are versions of the regular\vspace*{1pt} conditional probabilities
$\mathbf{P}((Y_k)_{k<0}\in\cdv \sigma(X_0)\vee\mathcal{F}^Y_+)$ and
$\mathbf{P}((Y_k)_{k<0}\in\cdv \mathcal{F}^Y_+)$, respectively.
Applying Lemma~\ref{lempolishbayes} to
$((Y_k)_{k\ge0},X_0,(Y_k)_{k<0})$, it suffices to show that
\[
\mathbf{P}^{z,y_0}\bigl((Y_k)_{k<0}\in\cdot\bigr)
\sim\mathbf{P}^{\varpi^+(y, \cdot)\otimes\delta_{y_0}}\bigl(
(Y_k)_{k<0}\in\cdot\bigr)
\]
for $\varpi\mathbf{P}^Y$-a.e. $(z,y)$. By Lemma~\ref{lemnondegrv},
we may apply Proposition~\ref{propcoupling} to the reverse-time model.
Therefore, it suffices to prove that
\[
\liminf_{n\to\infty}\bigl\|\mathbf{P}^{z,y_0}(X_{-n}\in\cdot)-
\mathbf{P}^{\varpi^+(y, \cdot)\otimes\delta_{y_0}}(
X_{-n}\in\cdot)\bigr\|_{\mathrm{TV}}=0
\]
for $\varpi\mathbf{P}^Y$-a.e. $(z,y)$. To this end, let us note that
\begin{eqnarray*}
&&\mathbf{E}\bigl(\bigl\|\mathbf{P}^{X_0,Y_0}(X_{-n}\in\cdot)-
\mathbf{P}^{\varpi^+(Y, \cdot)\otimes\delta_{Y_0}}(
X_{-n}\in\cdot)\bigr\|_{\mathrm{TV}}\bigr) \\
&&\qquad \le
\mathbf{E}\bigl(\|\mathbf{P}^{X_0,Y_0}(X_{-n}\in\cdot)-
\mathbf{P}(X_{-n}\in\cdot)\|_{\mathrm{TV}}\bigr) \\
&&\qquad\quad{} +
\mathbf{E}\bigl(\bigl\|\mathbf{P}(X_{-n}\in\cdv \mathcal{F}^Y_+)
-\mathbf{P}(X_{-n}\in\cdot)\bigr\|_{\mathrm{TV}}\bigr) \\
&&\qquad \le
2 \mathbf{E}\bigl(\|\mathbf{P}^{X_0,Y_0}(X_{-n}\in\cdot)-
\mathbf{P}(X_{-n}\in\cdot)\|_{\mathrm{TV}}\bigr).
\end{eqnarray*}
Thus the result follows by Assumption~\ref{asptergodicrv} and
Fatou's lemma.
\end{pf}

We now complete the proof of Theorem~\ref{thmexchg}.\vadjust{\goodbreak}
\begin{pf}
The first part of Theorem~\ref{thmexchg} follows immediately from
Proposition~\ref{propvw} and Lemma~\ref{lemmumuplus}. Now note that by
Lemma~\ref{lemnondegrv}, Assumptions
\mbox{\ref{asptergodic}--\ref{asptnondeg}} still hold if we
replace the model $(X_n,Y_n)_{n\in\mathbb{Z}}$ by the time-reversed
model $(X_{-n},Y_{-n})_{n\in\mathbb{Z}}$. Therefore, the second part of
Theorem~\ref{thmexchg} follows immediately from the first part by time
reversal.
\end{pf}

%s5 #&#
\section{\texorpdfstring{Proof of Theorem \protect\ref{thmfilter}}{Proof of Theorem 2.10}}
\label{secfiltstab}

The goal of this section is to prove Theorem~\ref{thmfilter}. We begin by
recalling some basic properties of the filter. Then, we prove
Theorem~\ref{thmfilter} first for a special case, then in
the general case by a~recursive argument.

%s5.1 #&#
\subsection{Preliminaries}

Recall that $\Pi_n^\mu$ is defined as a version of the regular conditional
probability $\mathbf{P}^\mu(X_n\in\cdv \mathcal{F}^Y_{[0,n]})$. Of
course, we are free to choose an arbitrary version of the filter, as the
statement of Theorem~\ref{thmfilter} does not depend on the choice of
version (this follows from Corollary~\ref{coryinfhabscont}).
Nonetheless, we will find it convenient in our proofs to work with
specific versions of these regular conditional probabilities, which we
define presently.

For notational simplicity, we introduce the following device: for
every probability measure $\rho$ on $E\times F$, we fix a probability
kernel $\rho_{\cdot}\dvtx F\times\mathcal{B}(E)\to[0,1]$ such that
$\rho_{Y_0}(A)=\mathbf{P}^\rho(X_0\in A|Y_0)$ for all $A\in\mathcal{B}(E)$
[i.e., $\rho_\cdot$ is a version of the regular conditional probability
$\mathbf{P}^\rho(X_0\in\cdv Y_0)$].\label{page30}
%Y_0)$}%
%
%le5.1 #&#
\begin{lem}
\label{lemrefprob}
Suppose that assumption~\ref{asptnondeg} holds. For every probability
measure $\mu$ on $E\times F$, we define a sequence of probability kernels
$\Pi^\mu_n\dvtx \Omega^Y\times\mathcal{B}(E)\to[0,1]$
($n\ge0$)
through the following recursion:
\begin{eqnarray*}
\Pi^\mu_n(y,A) &=& \frac{\int I_A(z)
g(z',y(n-1),z,y(n)) P_0(z',dz) \Pi^\mu_{n-1}(y,dz')}
{\int g(z',y(n-1),z,y(n)) P_0(z',dz) \Pi^\mu_{n-1}(y,dz')}, \\
\Pi^\mu_0(y,A) &=& \mu_{y(0)}(A).
\end{eqnarray*}
Then $\Pi^\mu_n$ is a version of the regular conditional probability
$\mathbf{P}^\mu(X_n\in\cdv \mathcal{F}^Y_{[0,n]})$ for every
$n\ge0$. Moreover, $\Pi^\mu_n(y, \cdot)\sim
\mathbf{P}^{\mu_{y(0)}\otimes\delta_{y(0)}}(X_n\in\cdot)$ for all
$y,n$.
\end{lem}
\begin{pf}
By construction, we have
\begin{eqnarray*}
&&\mathbf{P}^\mu(X_0\in dx_0,\ldots,X_n\in dx_n,Y_0\in
dy_0,\ldots,Y_n\in dy_n) \\
&&\qquad=\mu(E\times dy_0)
\mu_{y_0}(dx_0)
\prod_{i=0}^{n-1}g(x_i,y_i,x_{i+1},y_{i+1})
P_0(x_i,dx_{i+1})Q(y_i,dy_{i+1}).
\end{eqnarray*}
Therefore, the Bayes formula gives for any $A\in\mathcal{B}(E)$
\begin{eqnarray*}
&&\mathbf{P}^{\mu}\bigl(X_n\in A|\mathcal{F}^Y_{[0,n]}\bigr) \\
&&\qquad=\frac{\int I_A(x_n)\mu_{Y_0}(dx_0)
\prod_{i=0}^{n-1}g(x_i,Y_i,x_{i+1},Y_{i+1})P_0(x_i,dx_{i+1})}{
\int\mu_{Y_0}(dx_0)
\prod_{i=0}^{n-1}g(x_i,Y_i,x_{i+1},Y_{i+1})P_0(x_i,dx_{i+1})}.
\end{eqnarray*}
This clearly coincides with the recursive definition of $\Pi_n^\mu$.
Moreover, it follows directly that $\Pi_n^\mu(y,\cdot)\sim
\mu_{y(0)}P_0^n$ for all $y,n$. But note that
\[
\mathbf{P}^{\mu_{w}\otimes\delta_{w}}(X_n\in A) =
\int I_A(x_n)
\mu_{w}(dx_0)
f(w,x_0,\ldots,x_n)
\prod_{i=0}^{n-1}
P_0(x_i,dx_{i+1}),
\]
where we have defined
\[
f(y_0,x_0,\ldots,x_n) =
\int\prod_{i=0}^{n-1}g(x_i,y_i,x_{i+1},y_{i+1})
Q(y_i,dy_{i+1}).
\]
Therefore,
$\Pi^\mu_n(y, \cdot)\sim\mu_{y(0)}P_0^n\sim
\mathbf{P}^{\mu_{y(0)}\otimes\delta_{y(0)}}(X_n\in\cdot)$
for every $y,n$.
\end{pf}

Throughout the remainder of this section, the nonlinear filter $\Pi
_n^\mu$
will always be assumed to be chosen according to the particular version
defined in Lemma~\ref{lemrefprob}. This entails no loss of generality in
our final results.
%
%re5.2 #&#
\begin{rem}
From the recursive formula for $\Pi_n^\mu$, we can read off that
\[
\Pi_{n+m}^\mu(y,A) =
\Pi_m^{\Pi_n^\mu(y, \cdot)\otimes\delta_{y(n)}}(
\Theta^n y,A) \qquad\mbox{for all }n,m\ge0, y\in\Omega^Y,
A\in\mathcal{B}(E).
\]
This recursive property will play an important role in our proof.
One of the advantages of our specific choice of version of the
filter is that this property holds pathwise, so that we need not
worry about the joint measurability of $\Pi_n^{\mu}(y, \cdot)$ with
respect to $(y,\mu)$. Of course, our choice of version is not essential
and technicalities of this kind could certainly be resolved more
generally if one were so inclined.
\end{rem}

%s5.2 #&#
\subsection{The absolutely continuous case}

We begin by obtaining an explicit formula for the limit of
$\|\Pi_n^\mu-\Pi_n^\nu\|_{\mathrm{TV}}$ for absolutely continuous
measures $\mu\ll\nu$. This result will be applied recursively in the
proof of Theorem~\ref{thmfilter}.
%
%pr5.3 #&#
\begin{prop}
\label{proptvlimit}
For any probability measures $\mu,\nu$ on $E\times F$ with $\mu\ll\nu$
\begin{eqnarray*}
&&\mathbf{E}^\mu\Bigl[
{\limsup_{n\to\infty}}\|\Pi_n^\mu-\Pi_n^\nu\|_{\mathrm{TV}}
\Bigr] \\
&&\qquad=\mathbf{E}^\nu\biggl[\biggl|
\mathbf{E}^\nu\biggl(\frac{d\mu}{d\nu}(X_0,Y_0)\Big|
\bigcap_{n\ge0} \mathcal{F}^Y_+\vee\mathcal{F}^X_{[n,\infty[}
\biggr) -
\mathbf{E}^\nu\biggl(\frac{d\mu}{d\nu}(X_0,Y_0)\Big|
\mathcal{F}^Y_+\biggr)
\biggr|
\biggr].
\end{eqnarray*}
\end{prop}
\begin{pf}
As $d\mathbf{P}^\mu/d\mathbf{P}^\nu=(d\mu/d\nu)(X_0,Y_0)$ by the
Markov property, we have
\[
\mathbf{P}^\mu\bigl(X_n\in A|\mathcal{F}^Y_{[0,n]}\bigr) =
\frac{
\mathbf{E}^\nu(I_A(X_n)
\mathbf{E}^\nu(({d\mu}/{d\nu})(X_0,Y_0)|\sigma(X_n)\vee
\mathcal{F}^Y_{[0,n]}) |\mathcal{F}^Y_{[0,n]})
}{
\mathbf{E}^\nu(({d\mu}/{d\nu})(X_0,Y_0)|
\mathcal{F}^Y_{[0,n]})
}
\]
$\mathbf{P}^\mu$-a.s. by the Bayes formula. Therefore, we evidently have
\[
\frac{d\Pi_n^\mu}{d\Pi_n^\nu}(X_n)=
\frac{
\mathbf{E}^\nu(({d\mu}/{d\nu})(X_0,Y_0)|\sigma(X_n)\vee
\mathcal{F}^Y_{[0,n]})
}{
\mathbf{E}^\nu(({d\mu}/{d\nu})(X_0,Y_0)|
\mathcal{F}^Y_{[0,n]})
} , \qquad \mathbf{P}^\mu\mbox{-a.s.}
\]
In particular, we can write $\mathbf{P}^\mu$-a.s.
\[
\|\Pi_n^\mu-\Pi_n^\nu\|_{\mathrm{TV}} =
\int\biggl|\frac{d\Pi_n^\mu}{d\Pi_n^\nu}(x)-1\biggr|\Pi_n^\nu(dx) =
\frac{
\mathbf{E}^\nu(M_n|\mathcal{F}^Y_{[0,n]})
}{
\mathbf{E}^\nu(({d\mu}/{d\nu})(X_0,Y_0)
|\mathcal{F}^Y_{[0,n]})
},
\]
where we have defined
\[
M_n =
\biggl|
\mathbf{E}^\nu\biggl(\frac{d\mu}{d\nu}(X_0,Y_0)\Big|
\sigma(X_n)\vee\mathcal{F}^Y_{[0,n]}\biggr)
-
\mathbf{E}^\nu\biggl(\frac{d\mu}{d\nu}(X_0,Y_0)\Big|
\mathcal{F}^Y_{[0,n]}\biggr)
\biggr|.
\]
Thus it is easily seen that
\[
\mathbf{E}^\mu\Bigl[
{\limsup_{n\to\infty}}\|\Pi_n^\mu-\Pi_n^\nu\|_{\mathrm{TV}}
\Bigr]=
\mathbf{E}^\nu\Bigl[
\limsup_{n\to\infty}
\mathbf{E}^\nu\bigl(M_n|\mathcal{F}^Y_{[0,n]}\bigr)
\Bigr].
\]
Now note that, by the Markov property, $\mathcal{F}_{[n+1,\infty[}$
and $\sigma(X_0)\vee\mathcal{F}^Y_{[0,n-1]}$ are conditionally independent
given $\sigma(X_n,Y_n)$. Therefore,
\[
M_n = \biggl| \mathbf{E}^\nu\biggl(\frac{d\mu}{d\nu}(X_0,Y_0)\Big|
\mathcal{F}^Y_+\vee\mathcal{F}^X_{[n,\infty[}\biggr) -
\mathbf{E}^\nu\biggl(\frac{d\mu}{d\nu}(X_0,Y_0)\Big|
\mathcal{F}^Y_{[0,n]}\biggr) \biggr|.
\]
If $d\mu/d\nu$ were uniformly bounded, the result would follow
directly from the martingale convergence theorem and Hunt's lemma
(\cite{DM}, Theorem V.45).

In the case that $d\mu/d\nu$ is unbounded, define the truncated process
\[
M_n^k = \biggl| \mathbf{E}^\nu\biggl(\frac{d\mu}{d\nu}(X_0,Y_0)\wedge k\Big|
\mathcal{F}^Y_+\vee\mathcal{F}^X_{[n,\infty[}\biggr) -
\mathbf{E}^\nu\biggl(\frac{d\mu}{d\nu}(X_0,Y_0)\wedge k\Big|
\mathcal{F}^Y_{[0,n]}\biggr) \biggr|.
\]
By Hunt's lemma and dominated convergence,
\[
\lim_{k\to\infty}\lim_{n\to\infty}
\mathbf{E}^\nu\bigl(M_n^k|\mathcal{F}^Y_{[0,n]}\bigr) =
\mathbf{E}^\nu(M_\infty|\mathcal{F}^Y_+),\qquad
\mathbf{P}^\nu\mbox{-a.s.,}
\]
where $M_\infty=\lim_{n\to\infty}M_n$.
Therefore, we obtain $\mathbf{P}^\nu$-a.s.
\[
\limsup_{n\to\infty}
\mathbf{E}^\nu\bigl(M_n|\mathcal{F}^Y_{[0,n]}\bigr) =
\mathbf{E}^\nu(M_\infty|\mathcal{F}^Y_+) +
\limsup_{k\to\infty}
\limsup_{n\to\infty}
\mathbf{E}^\nu\bigl(M_n-M_n^k|\mathcal{F}^Y_{[0,n]}\bigr).
\]
It remains to note that the second term vanishes $\mathbf{P}^\nu$-a.s.,
\begin{eqnarray*}
&&\Bigl| \limsup_{k\to\infty} \limsup_{n\to\infty}
\mathbf{E}^\nu\bigl(M_n-M_n^k|\mathcal{F}^Y_{[0,n]}\bigr) \Bigr| \\
&&\qquad\le
2\limsup_{k\to\infty} \limsup_{n\to\infty} \mathbf{E}^\nu\biggl(
\frac{d\mu}{d\nu}(X_0,Y_0)- \frac{d\mu}{d\nu}(X_0,Y_0)\wedge k
\Big|\mathcal{F}^Y_{[0,n]}\biggr)=0.
\end{eqnarray*}
The proof is complete.
\end{pf}

%s5.3 #&#
\subsection{The general case}

In the special case where $\mu\ll\pi$, Theorem~\ref{thmfilter}
follows directly from Proposition~\ref{proptvlimit} and Theorem
\ref{thmexchg}. An additional step is needed, however, to prove
Theorem~\ref{thmexchg} in the general case.\vadjust{\goodbreak}
%
%le5.4 #&#
\begin{lem}
\label{lemfiltlebesgue}
Let $\mu,\rho$ be probability measures on $E$, and choose
$S\in\mathcal{B}(E)$ such that $\mu(S)>0$. Define
the probability measure $\nu=\mu( \cdot\cap S)/\mu(S)$. Then
\begin{eqnarray*}
\|\Pi^{\mu\otimes\delta_w}_n-
\Pi^{\rho\otimes\delta_w}_n\|_{\mathrm{TV}}
&\le&
2 \mathbf{P}^{\mu\otimes\delta_w}\bigl(X_0\notin S|
\mathcal{F}_{[0,n]}^Y\bigr)
\\
&&{}+
\mathbf{P}^{\mu\otimes\delta_w}\bigl(X_0\in S
|\mathcal{F}_{[0,n]}^Y\bigr)
\|\Pi^{\nu\otimes\delta_w}_n-
\Pi^{\rho\otimes\delta_w}_n\|_{\mathrm{TV}}
\end{eqnarray*}
$\mathbf{P}^{\mu\otimes\delta_w}$-a.s. for any $w\in F$.
\end{lem}
\begin{pf}
If $\mu(S)=1$, the proof is trivial. Otherwise, by the Bayes formula,
\[
\Pi^{\mu\otimes\delta_w}_n =
\mathbf{P}^{\mu\otimes\delta_w}\bigl(X_0\in S|\mathcal{F}_{[0,n]}^Y\bigr)
\Pi^{\nu\otimes\delta_w}_n +
\mathbf{P}^{\mu\otimes\delta_w}\bigl(X_0\notin S|\mathcal{F}_{[0,n]}^Y\bigr)
\Pi^{\nu^\perp\otimes\delta_w}_n
\]
$\mathbf{P}^{\mu\otimes\delta_w}$-a.s.,
where $\nu^\perp=\mu( \cdot\cap S^c)/\mu(S^c)$. But obviously
\[
\Pi^{\rho\otimes\delta_w}_n =
\mathbf{P}^{\mu\otimes\delta_w}\bigl(X_0\in S|\mathcal{F}_{[0,n]}^Y\bigr)
\Pi^{\rho\otimes\delta_w}_n +
\mathbf{P}^{\mu\otimes\delta_w}\bigl(X_0\notin S|\mathcal{F}_{[0,n]}^Y\bigr)
\Pi^{\rho\otimes\delta_w}_n
\]
$\mathbf{P}^{\mu\otimes\delta_w}$-a.s., so the result follows directly.
\end{pf}
%
%re5.5 #&#
\begin{rem}
Even though we have fixed a version of the filter
$\Pi_n^{\rho\otimes\delta_w}$, our results should ultimately not depend on
the choice of version. In this light, Lemma~\ref{lemfiltlebesgue} may appear somewhat
suspicious as the regular conditional probability
$\mathbf{P}^{\rho\otimes\delta_w}(X_n\in\cdv \mathcal{F}^Y_{[0,n]})$
is not $\mathbf{P}^{\mu\otimes\delta_w}$-a.s. uniquely defined. However,
there is no problem here, as the proof shows that the inequality in Lemma
\ref{lemfiltlebesgue} holds for \textit{any} choice of version, even though
different versions may be inequivalent. On the other hand, we will
ultimately apply this result only when $\nu\ll\rho$, in which case the
expression is in fact independent of the choice of version.
\end{rem}

The idea is now to apply the recursive property of the filter:
\[
\|\Pi_{m+n}^\mu-\Pi_{m+n}^\pi\|_{\mathrm{TV}} =
\|\Pi_n^{\Pi_m^\mu\otimes\delta_{Y_m}}(Y\circ\Theta^m, \cdot)-
\Pi_n^{\Pi_m^\pi\otimes\delta_{Y_m}}(Y\circ\Theta^m, \cdot)
\|_{\mathrm{TV}}
\]
for any $m\ge0$. As $\Pi_m^\mu\sim\mathbf{P}^\mu(X_m\in\cdv Y_0)$
and $\Pi_m^\pi\sim\mathbf{P}(X_m\in\cdv Y_0)$ by Lem\-ma~\ref{lemrefprob}, the assumption of Theorem~\ref{thmfilter} guarantees
that the singular part of~$\Pi_m^\mu$ with respect to $\Pi_m^\pi$
vanishes as $m\to\infty$. We can therefore use Lem\-ma~\ref{lemfiltlebesgue} to replace $\Pi_m^\mu$ by its absolutely continuous
part, so that we have reduced the limit as $n\to\infty$ to the special
case of Proposition~\ref{proptvlimit}. In order to apply Proposition
\ref{proptvlimit}, however, we will require one additional result.
%
%le5.6 #&#
\begin{lem}
\label{lemvwcondi}
Suppose that Assumptions~\ref{asptergodic}--\ref{asptnondeg} hold.
Then for any \mbox{$m\ge0$}
\[
\bigcap_{n\ge0} \mathcal{F}^Y_+\vee\mathcal{F}^X_{[n,\infty[}=
\mathcal{F}^Y_+,\qquad
\mathbf{P}^{\Pi_m^\pi(y, \cdot)\otimes\delta_{y(m)}}\mbox{-a.s.}\qquad
\mbox{for }\mathbf{P}^Y\mbox{-a.e. } y.
\]
\end{lem}
\begin{pf}
As in\vspace*{1pt} the proof of Proposition~\ref{propvw}, it suffices to establish
that $\mathcal{T}^X$ is
$\mathbf{P}^{\Pi_m^\pi(y, \cdot)\otimes\delta_{y(m)}}( \cdv
\mathcal{F}^Y_+)$-trivial
$\mathbf{P}^{\Pi_m^\pi(y, \cdot)\otimes\delta_{y(m)}}$-a.s. for
$\mathbf{P}^Y$-a.e. $y$. Note that
\[
\mathbf{P}^{\Pi_m^\pi(Y, \cdot)\otimes\delta_{Y_m}}(A) =
\mathbf{E}\bigl(\mathbf{P}^{X_m,Y_m}(A)|\mathcal{F}^Y_{[0,m]}\bigr) =
\mathbf{P}\bigl(A\circ\Theta^m|\mathcal{F}^Y_{[0,m]}\bigr)\vadjust{\goodbreak}
\]
for all $A\in\mathcal{F}_+$ by the Markov property.
Therefore,
\[
\mathbf{P}^{\Pi_m^\pi(Y\circ\Theta^{-m}, \cdot)
\otimes\delta_{Y_0}}(A) =
\mathbf{P}\bigl(A\circ\Theta^m|\mathcal{F}^Y_{[0,m]}\bigr)
\circ\Theta^{-m} =
\mathbf{P}\bigl(A|\mathcal{F}^Y_{[-m,0]}\bigr),
\]
where we used that $\mathbf{P}$ is stationary. It follows that
$\mathbf{P}^{\Pi_m^\pi(Y\circ\Theta^{-m}, \cdot)\otimes\delta
_{Y_0}}|_{\mathcal{F}_+}$
is a version of $\mathbf{P}((X_n,Y_n)_{n\ge0}\in\cdv
\mathcal{F}^Y_{[-m,0]})$. By Lemma~\ref{lemvWconditioning}
\[
\mathbf{P}\bigl((X_n)_{n\ge0}\in\cdv \mathcal{F}^Y_{[-m,\infty[}\bigr)
=
\mathbf{P}^{\Pi_m^\pi(Y\circ\Theta^{-m}, \cdot)\otimes\delta_{Y_0}}
\bigl((X_n)_{n\ge0}\in\cdv \mathcal{F}^Y_+\bigr),\qquad
\mathbf{P}\mbox{-a.s.}
\]
Thus it\vspace*{1pt} suffices to show that $\mathcal{T}^X$ is
$\mathbf{P}( \cdv \mathcal{F}^Y_{[-m,\infty[})$-trivial
$\mathbf{P}$-a.s., which is equivalent (by virtue of
Lemma~\ref{lemvWexchg} in Appendix~\ref{secapp}) to
\[
\bigcap_{n\ge0} \mathcal{F}^Y_{[-m,\infty[}\vee\mathcal{F}^X_{[n,\infty[}=
\mathcal{F}^Y_{[-m,\infty[} ,\qquad\mathbf{P}\mbox{-a.s.}
\]
But this follows directly from Theorem~\ref{thmexchg} and the
stationarity of $\mathbf{P}$.
\end{pf}

We can now complete the proof of Theorem~\ref{thmfilter}.
\begin{pf}
By the recursive property of the filter,
\begin{eqnarray*}
&&{\limsup_{k\to\infty}}
\|\Pi^\mu_{k}-\Pi^\pi_{k}\|_{\mathrm{TV}} \\
&&\qquad=\limsup_{k\to\infty}
\bigl\|\Pi_k^{\Pi^\mu_n(Y, \cdot)\otimes\delta_{Y_n}}(Y\circ\Theta^n,
\cdot)
- \Pi_k^{\Pi^\pi_n(Y, \cdot)\otimes\delta_{Y_n}}(Y\circ\Theta^n,
\cdot)
\bigr\|_{\mathrm{TV}}
\end{eqnarray*}
for all $n\ge0$. Therefore, we obtain $\mathbf{P}^\mu$-a.s.
\begin{eqnarray*}
&&\mathbf{E}^\mu\Bigl(
{\limsup_{k\to\infty}}\|\Pi^\mu_k-\Pi^\pi_k\|_{\mathrm{TV}}
|\mathcal{F}_{[0,n]}^Y\Bigr)
\\
&&\qquad=
\mathbf{E}^{\Pi^\mu_n(y, \cdot)\otimes\delta_{y(n)}}\Bigl(
\limsup_{k\to\infty}\bigl\|\Pi_k^{\Pi^\mu_n(y, \cdot)\otimes\delta_{y(n)}}
- \Pi_k^{\Pi^\pi_n(y, \cdot)\otimes\delta_{y(n)}}\bigr\|_{\mathrm{TV}}\Bigr)
\Big|_{y=Y},
\end{eqnarray*}
where we have used that $\Pi_n^\mu(Y, \cdot)$ and
$\Pi_n^\pi(Y, \cdot)$ are $\mathcal{F}^Y_{[0,n]}$-measurable.

To proceed, let us first recall that
\[
\Pi^\mu_n(y, \cdot)\sim
\mathbf{P}^{\mu_{y(0)}\otimes\delta_{y(0)}}(X_n\in\cdot)
\quad\mbox{and}\quad
\Pi^\pi_n(y, \cdot)\sim
\mathbf{P}^{\pi_{y(0)}\otimes\delta_{y(0)}}(X_n\in\cdot)
\]
for all $y,n$ by Lemma~\ref{lemrefprob}. Choose a set
$S_n\in\mathcal{B}(E\times
F)$ such that
\[
\mathbf{P}^{\mu_{w}\otimes\delta_{w}}\bigl(X_n\in\cdot\cap
S_n(w)\bigr)\ll\mathbf{P}^{\pi_{w}\otimes\delta_{w}}(X_n\in\cdot)
\]
and
\[
\mathbf{P}^{\pi_{w}\otimes\delta_{w}}\bigl(X_n\in S_n(w)\bigr)=1
\]
for all $w\in F$, where $I_{S_n(w)}(z) = I_{S_n}(z,w)$ (the existence of
such a set follows from the Lebesgue decomposition for kernels;~\cite{DM},
Section V.58). Define
\[
\Sigma_n(y, \cdot)=\Pi_n^\mu\bigl(y, \cdot\cap S_n(y(0))\bigr)/
\Pi_n^\mu(y,S_n(y(0))).
\]
Then clearly
$\Sigma_n(y, \cdot)\ll\Pi_n^\pi(y, \cdot)$ for all $y$, and
by Lemma~\ref{lemfiltlebesgue}
\begin{eqnarray*}
&& \mathbf{E}^{\Pi^\mu_n(y, \cdot)\otimes\delta_{y(n)}}\Bigl(
\limsup_{k\to\infty}\bigl\|\Pi_k^{\Pi^\mu_n(y, \cdot)\otimes\delta_{y(n)}}
- \Pi_k^{\Pi^\pi_n(y, \cdot)\otimes\delta_{y(n)}}\bigr\|_{\mathrm{TV}}\Bigr)
\\
&&\qquad
\le
2 \mathbf{P}^{\Pi^\mu_n(y, \cdot)\otimes\delta_{y(n)}}\bigl(
X_0\notin S_n(y(0))\bigr) \\
&&\qquad\quad{}  +
\mathbf{E}^{\Sigma_n(y, \cdot)\otimes\delta_{y(n)}}\Bigl(
\limsup_{k\to\infty}
\bigl\|\Pi_k^{\Sigma_n(y, \cdot)\otimes\delta_{y(n)}}
- \Pi_k^{\Pi^\pi_n(y, \cdot)\otimes\delta_{y(n)}}\bigr\|_{\mathrm{TV}}\Bigr).
\end{eqnarray*}
The last term vanishes for $\mathbf{P}^Y$-a.e. $y$ by Proposition
\ref{proptvlimit} and Lemma~\ref{lemvwcondi}, hence, for
$\mathbf{P}^\mu(Y\in\cdot)$-a.e. $y$ by Corollary
\ref{coryinfhabscont}. We have therefore shown that
\[
\mathbf{E}^\mu\Bigl(
{\limsup_{k\to\infty}}\|\Pi^\mu_k-\Pi^\pi_k\|_{\mathrm{TV}}
|\mathcal{F}_{[0,n]}^Y\Bigr) \le
2 \mathbf{P}^\mu\bigl(X_n\notin S_n(Y_0)|\mathcal{F}^Y_{[0,n]}\bigr),
\qquad\mathbf{P}^\mu\mbox{-a.s.}
\]
for every $n\ge0$. In particular, we have
\[
\mathbf{E}^\mu\Bigl(
{\limsup_{k\to\infty}}\|\Pi^\mu_k-\Pi^\pi_k\|_{\mathrm{TV}}\Bigr) \le
2 \mathbf{P}^\mu\bigl(X_n\notin S_n(Y_0)\bigr)\qquad
\mbox{for all }n\ge0.
\]
But as $\mathbf{P}(X_n\notin S_n(Y_0)|Y_0) =
\mathbf{P}^{\pi_{Y_0}\otimes\delta_{Y_0}}(X_n\notin S_n(Y_0))=0$,
we obtain
\begin{eqnarray*}
\mathbf{P}^\mu\bigl(X_n\notin S_n(Y_0)\bigr) &=&
\mathbf{E}^\mu\bigl(
\mathbf{P}^\mu\bigl(X_n\notin S_n(Y_0)|Y_0\bigr)-
\mathbf{P}\bigl(X_n\notin S_n(Y_0)|Y_0\bigr)\bigr) \\
&\le&
\mathbf{E}^\mu\bigl(
\bigl\|\mathbf{P}^\mu(X_n\in\cdv Y_0)-\mathbf{P}(X_n\in\cdv Y_0)
\bigr\|_{\mathrm{TV}}
\bigr)\xrightarrow{n\to\infty}0,
\end{eqnarray*}
where convergence follows as in the proof of Corollary
\ref{coryinfhabscont}. Therefore,
\[
{\limsup_{k\to\infty}}\|\Pi^\mu_k-\Pi^\pi_k\|_{\mathrm{TV}}=0,\qquad
\mathbf{P}^\mu\mbox{-a.s.},
\]
which completes the main part of the proof. To obtain $\mathbf
{P}$-a.s.
convergence (rather than $\mathbf{P}^\mu$-a.s. convergence) in the case
where $\mu(E\times\cdot)\sim\pi(E\times\cdot)$, it suffices to note
that in this case
$\mathbf{P}^\mu|_{\mathcal{F}^Y_+}\sim\mathbf{P}|_{\mathcal{F}^Y_+}$ by
Corollary~\ref{coryinfhabscont}.
\end{pf}

%s6 #&#
\section{\texorpdfstring{Proof of Theorem \protect\ref{thminvmeas}}{Proof of Theorem 2.12}}
\label{secinvmeas}

The goal of this section is to prove Theorem~\ref{thminvmeas}. We begin
by developing some details of the basic properties of
$(\Pi_n^\mu,Y_n)_{n\ge0}$ in Section~\ref{secintrofilt} under Assumption
\ref{asptnondeg}. We then complete the proof of Theorem
\ref{thminvmeas}.

%s6.1 #&#
\subsection{\texorpdfstring{Markov property of the pair $(\Pi_n^\mu,Y_n)_{n\ge0}$}
{Markov property of the pair (Pi n mu,Y n) n>=0}}
\label{secpairm}

Throughout this section, we assume that Assumption~\ref{asptnondeg} is in
force. We begin by defining a measurable map
$U\dvtx\mathcal{P}(E)\times F\times F\to\mathcal{P}(E)$ as follows:
\[
U(\nu,y_0,y_1)(A)
=\frac{\int I_A(z)g(z',y_0,z,y_1)P_0(z',dz)\nu(dz')}
{\int g(z',y_0,z,y_1)P_0(z',dz)\nu(dz')}.
\]
It follows immediately from Lemma~\ref{lemrefprob} that
$\Pi_n^\mu=U(\Pi_{n-1}^\mu,Y_{n-1},Y_n)$ $\mathbf{P}^\mu$-a.s. for every
$n\ge1$ and $\mu\in\mathcal{P}(E\times F)$.

Now define the transition kernel
$\Gamma\dvtx\mathcal{P}(E)\times F\times\mathcal{B}(\mathcal{P}(E)\times F)
\to[0,1]$ as
\[
\Gamma(\nu,y_0,A)=\int
I_A(U(\nu,y_0,y_1),y_1)P(z,y_0,dz',dy_1)\nu(dz).
\]
Then we have the following lemma.
%
%le6.1 #&#
\begin{lem}
\label{lemmarkf}
Suppose that Assumption~\ref{asptnondeg} holds. Then the
$(\mathcal{P}(E)\times F)$-valued process $(\Pi_n^\mu,Y_n)_{n\ge0}$ is
Markov under $\mathbf{P}^\mu$ with transition kernel $\Gamma$.
\end{lem}
\begin{pf}
It suffices to note that $(\Pi_n^\mu,Y_n)$
is $\mathcal{F}^Y_{[0,n]}$-measurable and
\begin{eqnarray*}
&&\mathbf{P}^\mu\bigl((\Pi^\mu_{n+1},Y_{n+1})\in
A|\mathcal{F}^Y_{[0,n]}\bigr)\\
&&\qquad=\mathbf{P}^\mu\bigl((U(\Pi^\mu_n,Y_n,Y_{n+1}),Y_{n+1})\in A|\mathcal{F}^Y_{[0,n]}\bigr)
\\
&&\qquad=\int I_A(U(\Pi^\mu_n,Y_n,w),w)P(z,Y_n,dz',dw)\Pi_n^\mu(dz)
=\Gamma(\Pi_n^\mu,Y_n,A)
\end{eqnarray*}
for every $A\in\mathcal{B}(\mathcal{P}(E)\times F)$.
\end{pf}

We can now establish some basic properties of $\Gamma$-invariant measures.
%
%le6.2 #&#
\begin{lem}
\label{lembarycenter}
Suppose that Assumption~\ref{asptnondeg} holds. Then for any
$\Gamma$-invar\-iant probability measure $\mathsf{m}$, the barycenter
$b\mathsf{m}$ is a $P$-invariant measure. Conversely, there is at least
one $\Gamma$-invariant measure with barycenter $\pi$.
\end{lem}
\begin{pf}
First, let $\mathsf{m}\in\mathcal{P}(\mathcal{P}(E)\times F)$ be a
$\Gamma$-invariant measure. Then
\begin{eqnarray*}
b\mathsf{m}(A\times B)
&=& \int\nu(A) I_B(w) \Gamma(\nu',w',d\nu,dw)
\mathsf{m}(d\nu',dw') \\
& = &
\int U(\nu',w',w)(A) I_B(w) P(z,w',d\tilde z,dw) \nu'(dz)
\mathsf{m}(d\nu',dw') \\
& = &
\int\frac{\int I_A(\tilde z)
g(z,w',\tilde z,w) P_0(z,d\tilde z) \nu'(dz)
}{\int g(z,w',\tilde z,w) P_0(z,d\tilde z) \nu'(dz)}
\\
&&\hspace*{9.5pt}{}\times
g(z,w',\tilde z,w) P_0(z,d\tilde z) \nu'(dz)
I_B(w) Q(w',dw) \mathsf{m}(d\nu',dw') \\
& = &
\int P(z,w',A\times B) \nu'(dz)
\mathsf{m}(d\nu',dw') \\
&=&
\int P(z,w',A\times B) b\mathsf{m}(dz,dw').
\end{eqnarray*}
Thus the barycenter $b\mathsf{m}$ is $P$-invariant.

To prove the converse, let $\Pi_n$ be a version of the regular conditional
probability $\mathbf{P}(X_n\in\cdv \mathcal{F}^Y_n)$, and let
$\Pi_{k,n}$ be a version of the regular conditional probability
$\mathbf{P}(X_n\in\cdv \mathcal{F}^Y_{[k,n]})$. Applying\vspace*{1pt} the Bayes
formula as in the proof of Lemma~\ref{lemrefprob}, we find that
$U(\Pi_{k,n},Y_{n},Y_{n+1})=\Pi_{k,n+1}$ $\mathbf{P}$-a.s. for every
$k\le n$. By the martingale convergence theorem, it follows directly that
\[
U(\Pi_{n},Y_{n},Y_{n+1})(A) =
\lim_{k\to-\infty}U(\Pi_{k,n},Y_{n},Y_{n+1})(A) = \Pi_{n+1}(A)
\]
$\mathbf{P}$-a.s. for every $A\in\mathcal{B}(E)$. As $\mathcal{B}(E)$ is
countably generated, a standard monotone class argument shows that
$U(\Pi_{n},Y_{n},Y_{n+1})=\Pi_{n+1}$ $\mathbf{P}$-a.s.
Therefore, the proof of Lemma~\ref{lemmarkf} shows that
$(\Pi_n,Y_n)_{n\in\mathbb{Z}}$ is Markov under $\mathbf{P}$ with
transition kernel $\Gamma$. But as $\mathbf{P}$ is stationary, the
process $(\Pi_n,Y_n)_{n\in\mathbb{Z}}$ is stationary
also. Therefore, the law of $(\Pi_0,Y_0)$ is a $\Gamma$-invariant
measure whose barycenter is obviously $\pi$.
\end{pf}

%s6.2 #&#
\subsection{\texorpdfstring{Uniqueness of the $\Gamma$-invariant measure}{Uniqueness of the Gamma-invariant measure}}

Given $\mathsf{m}\in\mathcal{P}(\mathcal{P}(E)\times F)$, define
the probability measure $\mathbf{P}_\mathsf{m}$ on the space
$\mathcal{P}(E)\times E^\mathbb{N} \times F^\mathbb{N}$ as
\begin{eqnarray*}
&&\mathbf{P}_{\mathsf{m}}\bigl((m_0,X_0,\ldots,X_n,Y_0,\ldots,Y_n)\in A\bigr)
\\
&&\qquad=\int I_A(\nu,x_0,\ldots,x_n,y_0,\ldots,y_n)
\nu(dx_0)P(x_0,y_0,dx_1,dy_1) \\
&&\qquad\quad\hspace*{9.4pt}{}\times\cdots\times
P(x_{n-1},y_{n-1},dx_n,dy_n)\mathsf{m}(d\nu,dy_0).
\end{eqnarray*}
We now choose regular versions of the following conditional
probabilites:
\begin{eqnarray*}
\sfPi
^{\min}_n &=& \mathbf{P}_\mathsf{m}\bigl(X_n\in\cdv
\mathcal{F}^Y_{[0,n]}\bigr),\\
\sfPi^\mathsf{m}_n &=& \mathbf{P}_\mathsf{m}\bigl(X_n\in\cdv
\sigma(m_0)\vee\mathcal{F}^Y_{[0,n]}\bigr),\\
\sfPi^{\max}_n &=& \mathbf{P}_\mathsf{m}\bigl(X_n\in\cdv
\sigma(m_0,X_0)\vee\mathcal{F}^Y_{[0,n]}\bigr).
\end{eqnarray*}
The following result is straightforward.
%
%le6.3 #&#
\begin{lem}
\label{lemkrandomiz}
The laws of $(\sfPi_n^{\min},Y_n)$ and $(\sfPi_n
^{\max},Y_n)$ under $\mathbf{P}_{\mathsf{m}}$ coincide with the laws of
$(\mathbf{P}^{b\mathsf{m}}(X_n\in\cdv \mathcal{F}^Y_{[0,n]}),
Y_n)$ and
$(\mathbf{P}^{b\mathsf{m}}(X_n\in\cdv \sigma(X_0)\vee\mathcal
{F}^Y_{[0,n]}),Y_n)$
under $\mathbf{P}^{b\mathsf{m}}$, respectively. Moreover, the process
$(\sfPi_n^{\mathsf{m}},Y_n)_{n\ge0}$ is Markov under~$\mathbf{P}_{\mathsf{m}}$ with transition kernel
$\Gamma$ and initial measure $\mathsf{m}$.
\end{lem}
\begin{pf}
By definition of the barycenter, the law of $(X_n,Y_n)_{n\ge0}$ under~$\mathbf{P}_{\mathsf{m}}$ coincides with the law of $(X_n,Y_n)_{n\ge0}$
under~$\mathbf{P}^{b\mathsf{m}}$. Moreover, it is easily seen that
$\sfPi^{\max}_n =\mathbf{P}_\mathsf{m}(X_n\in\cdv
\sigma(X_0)\vee\mathcal{F}^Y_{[0,n]})$ by the\vspace*{1pt} Markov property,
so $\sfPi^{\max}_n$ and $\sfPi^{\min}_n$ depend on
$(X_n,Y_n)_{n\ge0}$ only. This establishes the first part of the result.
The second part follows as in the proof of Lemma~\ref{lemmarkf}.
\end{pf}

We can now complete the proof of Theorem~\ref{thminvmeas}.
\begin{pf}
Throughout the proof, let $\mathsf{m}$ be a fixed $\Gamma$-invariant
probability measure with barycenter $\pi$.\vadjust{\goodbreak} We will show that, by virtue
of Theorem~\ref{thmexchg}, this invariant measure must necessarily
coincide with the invariant measure obtained in the proof of Lemma
\ref{lembarycenter}.

Let $p\in\mathbb{N}$, choose arbitrary bounded measurable functions
$f\dvtx E\to\mathbb{R}^p$ and $g\dvtx F\to\mathbb{R}$, and let
$\kappa\dvtx \mathbb{R}^{p+1}\to\mathbb{R}$ be a convex function.
Then $\kappa$ is necessarily continuous, so the function
$F\dvtx \mathcal{P}(E)\times F\to\mathbb{R}$ defined by
\[
F(\nu,w) = \kappa\biggl(g(w),
\int f(x) \nu(dx)\biggr)
\]
is bounded and measurable. By Jensen's inequality,
\[
\mathbf{E}_{\mathsf{m}}(F(\sfPi^{\min}_n,Y_n)) \le
\mathbf{E}_{\mathsf{m}}(F(\sfPi^\mathsf{m}_n,Y_n)) \le
\mathbf{E}_{\mathsf{m}}(F(\sfPi^{\max}_n,Y_n))
\]
for all $n\ge0$. Therefore, by Lemma~\ref{lemkrandomiz}
and the $\Gamma$-invariance of $\mathsf{m}$, we obtain
\begin{eqnarray*}
\mathbf{E}\bigl(
\kappa\bigl(g(Y_n),
\mathbf{E}\bigl(f(X_n)|\mathcal{F}^Y_{[0,n]}\bigr)
\bigr)\bigr) &\le&
\int F(\nu,w) \mathsf{m}(d\nu,dw) \\
&\le&
\mathbf{E}\bigl(
\kappa\bigl(g(Y_n),
\mathbf{E}\bigl(f(X_n)|\sigma(X_0)\vee\mathcal{F}^Y_{[0,n]}\bigr)
\bigr)\bigr).
\end{eqnarray*}
But using stationarity of $\mathbf{P}$ and the Markov property of
$(X_n,Y_n)_{n\in\mathbb{Z}}$,
\begin{eqnarray*}
\mathbf{E}\bigl(
\kappa\bigl(g(Y_n),
\mathbf{E}\bigl(f(X_n)|\mathcal{F}^Y_{[0,n]}\bigr)
\bigr)\bigr) &=&
\mathbf{E}\bigl(
\kappa\bigl(g(Y_0),
\mathbf{E}\bigl(f(X_0)|\mathcal{F}^Y_{[-n,0]}\bigr)
\bigr)\bigr),\\
\mathbf{E}\bigl(
\kappa\bigl(g(Y_n),
\mathbf{E}\bigl(f(X_n)|\sigma(X_0)\vee\mathcal{F}^Y_{[0,n]}\bigr)
\bigr)\bigr) &=&
\mathbf{E}\bigl(
\kappa\bigl(g(Y_0),
\mathbf{E}\bigl(f(X_0)|\mathcal{F}_0^Y\vee\mathcal{F}_{-n}^X\bigr)
\bigr)\bigr)
\end{eqnarray*}
for all $n\ge0$. Thus martingale convergence
and Theorem~\ref{thmexchg} yield
\[
\int F(\nu,w)\mathsf{m}(d\nu,dw) =
\mathbf{E}(
\kappa(g(Y_0),
\mathbf{E}(f(X_0)|\mathcal{F}^Y_{0})
)) =
\int F(\nu,w)\mathsf{m}^0(d\nu,dw),
\]
where $\mathsf{m}^0$ denotes the distinguished $\Gamma$-invariant measure
obtained in the proof of Lemma~\ref{lembarycenter}. But a standard
approximation argument shows that class of functions of the form
$F(\nu,w)=\kappa(g(w),\int f(x)\nu(dx))$ is measure-determining (see,
e.g., the proof of Proposition A.7~\cite{vH11}), so we can conclude that
$\mathsf{m}=\mathsf{m}^0$. Thus we have shown that any $\Gamma$-invariant
probability measure with barycenter $\pi$ must coincide with
$\mathsf{m}^0$, which establishes uniqueness.

To complete the proof, it remains to consider the case when $P$ has unique
invariant probability measure (i.e., $\pi$ is the only $P$-invariant
probability measure). As the barycenter of any $\Gamma$-invariant
probability measure must be $P$-invariant, this implies that any
$\Gamma$-invariant measure must have barycenter~$\pi$. Therefore, in this
case, $\Gamma$ has a unique invariant probability measure.
\end{pf}

%s7 #&#
\section{\texorpdfstring{Proof of Theorem \protect\ref{thmlambdainv}}{Proof of Theorem 2.13}}
\label{seclambdainv}

The goal of this section is to prove Theorem~\ref{thmlambdainv}. We
begin by developing some details of the basic properties of
$(\Pi_n^\mu,X_n,Y_n)_{n\ge0}$ in Section~\ref{secintrofilt} under
Assumption~\ref{asptnondeg}. We then complete the proof of
Theorem~\ref{thmlambdainv}.

%s7.1 #&#
\subsection{\texorpdfstring{Markov property of the triple $(\Pi_n^\mu,X_n,Y_n)_{n\ge0}$}
{Markov property of the triple (Pi n mu,X n,Y n) n>=0}}
\label{sectriplem}

In this section we use the notation of Section
\ref{secpairm}, and we again assume that Assumption
\ref{asptnondeg}\vadjust{\goodbreak} is in force. Define the transition kernel
$\Lambda\dvtx\mathcal{P}(E)\times E\times F\times
\mathcal{B}(\mathcal{P}(E)\times E\times F)\to[0,1]$ as
\[
\Lambda(\nu,x_0,y_0,A)=\int
I_A(U(\nu,y_0,y_1),x_1,y_1)P(x_0,y_0,dx_1,dy_1).
\]
Then we have the following lemma.
%
%le7.1 #&#
\begin{lem}
\label{lemmarksgf}
Suppose that Assumption~\ref{asptnondeg} holds. Then
$(\Pi_n^\mu,X_n,Y_n)_{n\ge0}$ is a $(\mathcal{P}(E)\times E\times
F)$-valued Markov chain under $\mathbf{P}^\mu$ with transition kernel~$\Lambda$.
\end{lem}

\begin{pf}
It suffices to note that $(\Pi_n^\mu,X_n,Y_n)$
is $\mathcal{F}_{[0,n]}$-measurable and
\begin{eqnarray*}
&&\mathbf{P}^\mu\bigl((\Pi^\mu_{n+1},X_{n+1},Y_{n+1})\in
A|\mathcal{F}_{[0,n]}\bigr)\\
&&\qquad =\int I_A(U(\Pi^\mu_n,Y_n,w),z,w)P(X_n,Y_n,dz,dw)\\
&&\qquad=\Lambda(\Pi_n^\mu,X_n,Y_n,A)
\end{eqnarray*}
for every $A\in\mathcal{B}(\mathcal{P}(E)\times E\times F)$.
\end{pf}

For any probability measure
$\mathsf{M}\in\mathcal{P}(\mathcal{P}(E)\times E\times F)$, we
define probability measures $m\mathsf{M}\in\mathcal{P}(E\times F)$
and $\gamma\mathsf{M}\in\mathcal{P}(\mathcal{P}(E)\times F)$ as follows:
\[
m\mathsf{M}(A\times B) =
\mathsf{M}\bigl(\mathcal{P}(E)\times A\times B\bigr),\qquad
\gamma\mathsf{M}(C\times B) =
\mathsf{M}(C\times E\times B).
\]
Moreover, we define the class
\begin{eqnarray*}
\mathfrak{M}&=&
\biggl\{\mathsf{M}\in\mathcal{P}\bigl(\mathcal{P}(E)\times E\times F\bigr)\dvtx
\forall
A\in\mathcal{B}(\mathcal{P}(E)),
B\in\mathcal{B}(E),
C\in\mathcal{B}(F), \\
&&\hspace*{55pt}\mathsf{M}(A\times B\times C) =
\int\nu(B) I_{A\times C}(\nu,w) \mathsf{M}(d\nu,dz,dw)
\biggr\}.
\end{eqnarray*}
We can now establish some basic properties of $\Lambda$-invariant
measures.
%
%le7.2 #&#
\begin{lem}
\label{lemmarginal}
Suppose that Assumption~\ref{asptnondeg} holds. Then for any
$\Lambda$-invar\-iant probability measure $\mathsf{M}$, the marginal
$m\mathsf{M}$ is a $P$-invariant measure. If, in addition,
$\mathsf{M}\in\mathfrak{M}$, then $\gamma\mathsf{M}$
is a $\Gamma$-invariant measure with barycenter $m\mathsf{M}$.
Conversely, there is at least one $\Lambda$-invariant
$\mathsf{M}\in\mathfrak{M}$ with marginal $\pi$.
\end{lem}
\begin{pf}
Let $\mathsf{M}\in\mathcal{P}(\mathcal{P}(E)\times E\times F)$ be a
$\Lambda$-invariant probability measure. It is trivial that
$m\mathsf{M}$ is $P$-invariant. Now suppose that also
$\mathsf{M}\in\mathfrak{M}$. Then
\begin{eqnarray*}
\gamma\mathsf{M}(A) &=&
\int I_A(\nu,w) \mathsf{M}(d\nu,dz,dw) \\
&=&
\int
I_A(\nu',w') \Lambda(\nu,z,w,d\nu',dz',dw') \mathsf{M}(d\nu,dz,dw)
\\
&=&
\int I_A(U(\nu,w,w'),w') P(z,w,dz',dw') \mathsf{M}(d\nu,dz,dw) \\
&=&
\int
I_A(U(\nu,w,w'),w') P(z,w,dz',dw') \nu(dz) \gamma\mathsf{M}(d\nu,dw)
\\
&=&
\int I_A(\nu',w') \Gamma(\nu,w,d\nu',dw') \gamma\mathsf{M}(d\nu,dw),
\end{eqnarray*}
where we have used that $\mathsf{M}\in\mathfrak{M}$ in the penultimate
equality. Thus $\gamma\mathsf{M}$ is a~$\Gamma$-invariant measure.
Moreover, it follows from the definition of $\mathfrak{M}$ that
\[
\int\nu(B) I_C(w) \gamma\mathsf{M}(d\nu,dw) =
\mathsf{M}\bigl(\mathcal{P}(E)\times B\times C\bigr) =
m\mathsf{M}(B\times C),
\]
so $m\mathsf{M}$ is the barycenter of $\gamma\mathsf{M}$.
Finally, let $\Pi_0$ be a version of the regular conditional
probability $\mathbf{P}(X_0\in\cdv \mathcal{F}^Y_0)$. Then as
in the proof of Lemma~\ref{lembarycenter}, the law of
$(\Pi_0,X_0,Y_0)$ is a $\Lambda$-invariant measure in $\mathfrak{M}$ with
marginal $\pi$.
\end{pf}

%s7.2 #&#
\subsection{\texorpdfstring{Uniqueness of the $\Lambda$-invariant measure}{Uniqueness of the Lambda-invariant measure}}

The first part of the proof of Theorem~\ref{thmlambdainv} follows
easily from Theorem~\ref{thminvmeas} and Lemma~\ref{lemmarginal}.
%
%le7.3 #&#
\begin{lem}
\label{lemlambdafirstpart}
Suppose that Assumptions~\ref{asptergodic}--\ref{asptnondeg} hold.
Then there is\break a~unique $\Lambda$-invariant probability measure with
marginal $\pi$ in the class $\mathfrak{M}$.
\end{lem}
\begin{pf}
Lemma~\ref{lemmarginal} guarantees the existence of a $\Lambda$-invariant
measure in $\mathfrak{M}$ with marginal $\pi$. To prove uniqueness,
note that every probability measure $\mathsf{M}\in\mathfrak{M}$ is
uniquely determined by $\gamma\mathsf{M}$ as
\[
\mathsf{M}(A\times B\times C) =
\int\nu(B) I_{A\times C}(\nu,w) \gamma\mathsf{M}(d\nu,dw).
\]
Therefore, by Lemma~\ref{lemmarginal}, if there were to exist two
distinct $\Lambda$-invariant measures in $\mathfrak{M}$ with marginal
$\pi$, then there must exist two distinct $\Gamma$-invariant measures with
barycenter $\pi$, in contradiction with Theorem~\ref{thminvmeas}.
\end{pf}

The second part of the proof of Theorem~\ref{thmlambdainv} relies on
Theorem~\ref{thmfilter} instead of Theorem~\ref{thminvmeas}. To prepare
for the proof, we begin by showing that the strengthened variant of
Assumption~\ref{asptergodic} in Theorem~\ref{thmlambdainv} is
equivalent to the requirement that the assumption of Theorem
\ref{thmfilter} holds universally.
%
%le7.4 #&#
\begin{lem}
\label{lemuniversalerg}
The following are equivalent:
\begin{longlist}[(2)]
\item[(1)] For every probability measure $\mu$ on $E\times F$ such that
$\mu(E\times\cdot)=\pi(E\times\cdot)$
\[
\int
\|\mathbf{P}^{z,w}(X_n\in\cdot)-\mathbf{P}(X_n\in\cdot)
\|_{\mathrm{TV}} \mu(dz,dw)\xrightarrow{n\to\infty}0.
\]
\item[(2)] For every probability measure $\mu$ on $E\,{\times}\,F$ such that
\mbox{$\mu(E\,{\times}\,\cdot)\,{\ll}\,\pi(E\,{\times}\,\cdot)$}
\[
\mathbf{E}^\mu\bigl(
\bigl\|\mathbf{P}^\mu(X_n\in\cdv Y_0)-\mathbf{P}(X_n\in\cdot)
\bigr\|_{\mathrm{TV}}\bigr)\xrightarrow{n\to\infty}0.
\]
\end{longlist}
\end{lem}
\begin{pf}
$(1)\Rightarrow(2)$. Let $\mu$ be any probability measure on $E\times F$
with $\mu(E\times\cdot)\ll\pi(E\times\cdot)$, let $\mu_w(dz)$
be a version of the regular conditional probability
$\mathbf{P}^\mu(X_0\in\cdv Y_0)$ and define $\mu'(dz,dw)=
\mu_w(dz)\pi(E\times dw)$. Then $\mu'(E\times\cdot)=
\pi(E\times\cdot)$, so the first statement of the lemma implies that
we have
\[
\|\mathbf{P}^{X_0,Y_0}(X_n\in\cdot)-\mathbf{P}(X_n\in\cdot)
\|_{\mathrm{TV}}\xrightarrow{n\to\infty}0\qquad
\mbox{in }\mathbf{P}^{\mu'}\mbox{-probability}.
\]
But $\mu\ll\mu'$ by construction, so the convergence also holds in
$\mathbf{P}^\mu$-probability. Therefore, we obtain by dominated
convergence
\begin{eqnarray*}
&&\mathbf{E}^\mu\bigl(
\bigl\|\mathbf{P}^\mu(X_n\in\cdv Y_0)-\mathbf{P}(X_n\in\cdot)
\bigr\|_{\mathrm{TV}}\bigr) \\
&&\qquad
\le
\mathbf{E}^\mu\bigl(
\|\mathbf{P}^{X_0,Y_0}(X_n\in\cdot)-\mathbf{P}(X_n\in\cdot)
\|_{\mathrm{TV}}\bigr)\xrightarrow{n\to\infty}0.
\end{eqnarray*}
Thus the second statement of the lemma follows.

$(2)\Rightarrow(1)$.
Let $\mu$ be any probability measure on $E\times F$ such that
$\mu(E\times\cdot)=\pi(E\times\cdot)$ and let $\mu_w(dz)$
be a version of the regular conditional probability
$\mathbf{P}^\mu(X_0\in\cdv Y_0)$. By~\cite{Kal02}, Lemma 3.22,
there is a measurable function $\iota\dvtx F\times[0,1]\to E$ such that
$\int f(z) \mu_w(dz) = \int_0^1 f(\iota(w,x)) \,dx$ for all $w$.
Applying the second statement of the lemma to
$\mu^x(dz,dw) = \delta_{\iota(w,x)}(dz)\mu(E\times dw)=
\delta_{\iota(w,x)}(dz)\pi(E\times dw)$ gives
\[
\int
\bigl\|\mathbf{P}^{\iota(w,x),w}(X_n\in\cdot)-\mathbf{P}(X_n\in\cdot)
\bigr\|_{\mathrm{TV}} \mu(E\times dw)\xrightarrow{n\to\infty}0
\qquad\mbox{for all }x\in[0,1].
\]
Thus the first statement of the lemma follows by integrating with respect
to $\int_0^1\cdot \,dx$ and applying the dominated convergence theorem.
\end{pf}

Let us note that only the first half of this result is needed in what
follows. However, the equivalence of the two assumptions shows that we
have not unnecessarily strengthened the assumptions of Theorem
\ref{thmlambdainv}.

For the proof of Theorem~\ref{thmlambdainv}, we require another lemma.
%
%le7.5 #&#
\begin{lem}
\label{lemcvginclsig}
Suppose that Assumptions~\ref{asptergodic}--\ref{asptnondeg} are in
force and that
\[
\mathbf{E}^\mu\bigl(
\bigl\|\mathbf{P}^\mu(X_n\in\cdv Y_0)-\mathbf{P}(X_n\in\cdot)
\bigr\|_{\mathrm{TV}}\bigr)\xrightarrow{n\to\infty}0
\]
for every probability measure $\mu$ on $E\times F$ with
$\mu(E\times\cdot)\ll\pi(E\times\cdot)$. Then
\[
\int\mathbf{E}^{z,w}\bigl(\bigl\|\Pi_n^{m(z,w)\otimes\delta_w}-
\Pi_n^\pi\bigr\|_{\mathrm{TV}}\bigr) \pi(dz,dw)\xrightarrow{n\to\infty}0
\]
for any measurable function $m\dvtx E\times F\to\mathcal{P}(E)$.
\end{lem}
\begin{pf}
By Proposition~\ref{propproductpi} and the Bayes formula, there is a
strictly positive measurable function $h\dvtx E\times F\to\mathbb{R}_+$ such
that the probability kernel
\[
\pi^X(z,A) = \frac{\int I_A(w) h(z,w) \pi(E\times dw)}
{\int h(z,w) \pi(E\times dw)}\qquad
\mbox{for all }z\in E, A\in\mathcal{B}(F)
\]
is a version of the regular conditional probability
$\mathbf{P}(Y_0\in\cdv X_0)$. In particular, $\pi^X(z, \cdot)
\sim\pi(E\times\cdot)$ for all $z\in E$, so by our assumptions and
Corollary~\ref{coryinfhabscont} we obtain
$\mathbf{P}^{\delta_z\otimes\pi^X(z, \cdot)}|_{\mathcal{F}^Y_+}\sim
\mathbf{P}|_{\mathcal{F}^Y_+}$ for all $z\in E$.

Fix a measurable function $m\dvtx E\times F\to\mathcal{P}(E)$.
For every $z\in E$, define $\mu^{z}(dz',dw)=m(z,w)(dz')\pi(E\times dw)$.
Then by Theorem~\ref{thmfilter}, we have
\[
\|\Pi^{\mu^z}_n-\Pi^\pi_n\|_{\mathrm{TV}}
\xrightarrow{n\to\infty}0,\qquad
\mathbf{P}\mbox{-a.s.}
\]
for all $z\in E$. Thus by
$\mathbf{P}^{\delta_z\otimes\pi^X(z, \cdot)}|_{\mathcal{F}^Y_+}\sim
\mathbf{P}|_{\mathcal{F}^Y_+}$ and dominated convergence,
\[
\int\mathbf{E}^{z,w}(\|\Pi^{\mu^z}_n-\Pi^\pi_n\|_{\mathrm{TV}})
\pi^X(z,dw)
\xrightarrow{n\to\infty}0
\]
for all $z\in E$. But by Lemma~\ref{lemrefprob} we have
$\Pi^{\mu^z}_n=\Pi^{m(z,w)\otimes\delta_w}_n$ $\mathbf{P}^{z,w}$-a.s.
for all $n\ge0$. Integrating with respect to $\pi(dz\times F)$ and
applying the dominated convergence theorem completes the proof.
\end{pf}

We now proceed to the proof of Theorem~\ref{thmlambdainv}. Let $Z$ be
any Polish space endowed with the complete metric $d_Z$. Recall that
the space $\mathcal{P}(Z)$ is Polish when endowed with the metric (cf.
\cite{Dud02}, Theorem 11.3.3 and Corollary~11.5.5)
\begin{eqnarray*}
d_{\mathcal{P}(Z)}(\nu,\nu')&=& \sup\biggl\{
\biggl|\int f(z) \nu(dz)-\int f(z) \nu'(dz)\biggr|\dvtx
\sup_{x\in Z}|f(x)|\le1,\\
&&\hspace*{124pt}
\sup_{x,y\in Z}\frac{|f(x)-f(y)|}{d_Z(x,y)}\le1
\biggr\}.
\end{eqnarray*}
In particular, the complete metric
\[
D((\nu,z,w),(\nu',z',w'))=
d_{\mathcal{P}(E)}(\nu,\nu')+d_{E}(z,z')+d_F(w,w')
\]
metrizes the topology of $\mathcal{P}(E)\times E\times F$.
\begin{pf*}{Proof of Theorem~\ref{thmlambdainv}}
The first part of the theorem was established in Lemma
\ref{lemlambdafirstpart}. For the remainder of the proof, let us assume
that one of the equivalent assumptions in Lemma~\ref{lemuniversalerg} is
in force. We will show that any two $\Lambda$-invariant probability
measures with marginal $\pi$ must coincide.\looseness=1

To this end, let $\mathsf{M}$ and $\mathsf{M}'$ be two $\Lambda$-invariant
probability measures with marginal~$\pi$. By~\cite{Kal02}, Lemma 3.22,
there exist measurable functions\vadjust{\goodbreak} $m\dvtx \break E\times F\times
[0,1]\to\mathcal{P}(E)$ and $m'\dvtx  E\times F\times[0,1]\to\mathcal{P}(E)$
such that
\begin{eqnarray*}
\int f(\nu,z,w) \mathsf{M}(d\nu,dz,dw)&=&
\int_0^1\int f(m(z,w,x),z,w) \pi(dz,dw) \,dx, \\
\int f(\nu,z,w) \mathsf{M}'(d\nu,dz,dw)&=&
\int_0^1\int f(m'(z,w,x),z,w) \pi(dz,dw) \,dx
\end{eqnarray*}
for every bounded measurable function $f\dvtx \mathcal{P}(E)\times
E\times F\to\mathbb{R}$. Moreover, note that by the definition of
$\Lambda$ and Lemma~\ref{lemrefprob}
\[
\int f(\nu',z',w') \Lambda^n(\nu,z,w,d\nu',dz',dw') =
\mathbf{E}^{z,w}(f(\Pi_n^{\nu\otimes\delta_w},X_n,Y_n)).
\]
Let us now fix a bounded function $f$ such that
\[
|f(\nu,z,w)-f(\nu',z',w')|\le
D((\nu,z,w),(\nu',z',w'))
\]
for all $\nu,\nu'\in\mathcal{P}(E)$,
$z,z'\in E$, $w,w'\in F$. We can now estimate
\begin{eqnarray*}
&&\biggl|\int f(\nu,z,w) \mathsf{M}(d\nu,dz,dw)-
\int f(\nu,z,w) \mathsf{M}'(d\nu,dz,dw)\biggr|
\\
&&\qquad \le
\int_0^1\int
\mathbf{E}^{z,w}\bigl(d_{\mathcal{P}(E)}\bigl(
\Pi_n^{m(z,w,x)\otimes\delta_w},\Pi_n^{m'(z,w,x)\otimes\delta_w}\bigr)\bigr)
\pi(dz,dw) \,dx
\\
&&\qquad \le
\int_0^1\int
\mathbf{E}^{z,w}\bigl(\bigl\|\Pi_n^{m(z,w,x)\otimes\delta_w}-
\Pi_n^{m'(z,w,x)\otimes\delta_w}\bigr\|_{\mathrm{TV}}\bigr)
\pi(dz,dw) \,dx
\end{eqnarray*}
for every $n\ge0$, where we used that $\mathsf{M}\Lambda^n=\mathsf{M}$
and $\mathsf{M}'\Lambda^n=\mathsf{M}'$. By the triangle inequality,
Lemma~\ref{lemcvginclsig} and the dominated convergence theorem, the
right-hand side of this inequality converges to zero as $n\to\infty$.
Therefore, we have shown that
\[
\biggl|\int f(\nu,z,w) \mathsf{M}(d\nu,dz,dw)-
\int f(\nu,z,w) \mathsf{M}'(d\nu,dz,dw)\biggr|=0
\]
for all bounded functions $f$ that are $1$-Lipschitz for the metric $D$.
In other words, $d_{\mathcal{P}(\mathcal{P}(E)\times E\times
F)}(\mathsf{M},\mathsf{M}')=0$, so $\mathsf{M}=\mathsf{M}'$.
Thus we have shown that all $\Lambda$-invariant
probability measures with marginal $\pi$ must coincide,
establishing uniqueness.

To complete the proof, it remains to consider the case when $P$ has unique
invariant probability measure (i.e., $\pi$ is the only $P$-invariant
probability measure). As the marginal of any $\Lambda$-invariant
probability measure must be $P$-invariant, this implies that any
$\Lambda$-invariant measure must have marginal~$\pi$. Therefore, in
this case, $\Lambda$ has a unique invariant probability measure.
\end{pf*}
%
%re7.6 #&#
\begin{rem}
It is instructive to note that Assumptions
\ref{asptergodic}--\ref{asptnondeg} are not sufficient to ensure
uniqueness of the $\Lambda$-invariant probability measure even in
the\vadjust{\goodbreak}
case that $P$ has a unique invariant probability measure.
Let us briefly sketch a counterexample. Let $E=\mathbb{R}\times\{0,1\}$
and $F=\mathbb{R}$, and consider the filtering model
\[
X_n^1 = 2X_{n-1}^1X_{n-1}^2 + \xi_n,\qquad
X_n^2 = X_{n-1}^2,\qquad
Y_n = X_n^1 + \eta_n,
\]
where $(\xi_n)_{n\ge0}$, $(\eta_n)_{n\ge0}$ are i.i.d.
$N(0,1)$-distributed random variables. It is clear that the
corresponding transition kernel $P$ has a unique invariant probability
measure $\pi$ [with $\pi( \cdot\times F)=N(0,1)\otimes\delta_0$] and
that Assumptions~\ref{asptergodic}--\ref{asptnondeg} hold.

Now let $\mu= \delta_0\otimes\delta_1\otimes N(0,1)$. Then
$\Pi_n^\mu=N(m_n,\sigma_n^2)\otimes\delta_1$, where $m_n$
and~$\sigma_n^2$ can be computed recursively using the Kalman filtering
equations corresponding to the model $X_n=2X_{n-1}+\xi_n$, $Y_n = X_n +
\eta_n$. It is easily verified by inspection of the Kalman filtering
equations that the law of $(\Pi_n^\mu,X_n,Y_n)$ converges weakly as
$n\to\infty$ under the stationary measure $\mathbf{P}$. The limiting
law is therefore a $\Lambda$-invariant probability measure that is
supported on $\mathcal{P}(\mathbb{R}\times\{1\})\times E\times F$. On
the other hand, the $\Lambda$-invariant measure defined in the proof of
Lemma~\ref{lemmarginal} is clearly supported on
$\mathcal{P}(\mathbb{R}\times\{0\})\times E\times F$. Therefore,
$\Lambda$ has distinct invariant measures.

This example illustrates that the stronger assumption of Theorem
\ref{thmlambdainv} is indeed required to establish uniqueness
of the $\Lambda$-invariant measure in
the class of all probability measures. Of course, the first part
of Theorem~\ref{thmlambdainv} is not contradicted as the additional
$\Lambda$-invariant measure obtained in this example is not in
$\mathfrak{M}$.
\end{rem}

\begin{appendix}
%s8 #&#
\section{Auxiliary results}
\label{secapp}

The goal of this Appendix is to collect for easy reference a few
auxiliary results that are used throughout the paper.

The following result on the existence of invariant sets for stationary
Markov chains is given in~\cite{vH10}, Lemma 2.6. The construction of
the set~$H$ follows closely along the lines of
\cite{OS70}, pages 1636 and 1637, so the proof is omitted.
%
%le8.1 #&#
\begin{lem}
\label{lemmcreinvset}
Let $\mathbf{P}^z$ be the law of a Markov process $(Z_k)_{k\ge0}$ given
\mbox{$Z_0=z$}, and let $\nu$ be a stationary probability for this Markov
process. Then for any set $\tilde H$ of $\nu$-full measure, there is a
subset $H\subset\tilde H$ of $\nu$-full measure such that
\[
\mathbf{P}^z(Z_n\in H\mbox{ for all }n\ge0)=1\qquad
\mbox{for all }z\in H.
\]
\end{lem}

The following elementary can be found in~\cite{vH10}, Lemma 3.6.
%
%le8.2 #&#
\begin{lem}
\label{lempolishbayes}
$\!\!\!$Let $G_1$, $G_2$ and $K$ be Polish spaces
and set \mbox{$\Omega\,{=}\,G_1\,{\times}\,G_2\,{\times}\,K$}.
We consider a probability measure
$\mathbf{P}$ on $(\Omega,\mathcal{B}(\Omega))$. Denote by
$\gamma_1\dvtx\Omega\to G_1$, $\gamma_2\dvtx\Omega\to G_2$, and
$\kappa\dvtx\Omega\to K$ the coordinate projections, and let
$\mathcal{G}_1$, $\mathcal{G}_2$, and $\mathcal{K}$ be the
$\sigma$-fields generated by $\gamma_1$, $\gamma_2$ and $\kappa$,
respectively. Choose fixed versions of the following regular
conditional probabilities:
\begin{eqnarray*}
\Xi^K_1(g_1,\cdot) &=& \mathbf{P}(\kappa\in\cdv \mathcal{G}_1)(g_1),
\qquad \Xi_{12}^K(g_1,g_2,\cdot) =
\mathbf{P}(\kappa\in\cdv \mathcal{G}_1\vee\mathcal{G}_2)(g_1,g_2),\\
\Xi^2_1(g_1,\cdot) &=&
\mathbf{P}(\gamma_2\in\cdv \mathcal{G}_1)(g_1),\qquad
\Xi_{1K}^2(g_1,k,\cdot) =
\mathbf{P}(\gamma_2\in\cdv \mathcal{G}_1\vee\mathcal{K})(g_1,k),
\end{eqnarray*}
where $g_1\in G_1$, $g_2\in G_2$, $k\in K$. Suppose that there
exists a nonnegative measurable function $h\dvtx G_1\times G_2\times
K\to[0,\infty[$ and a set $H\subset G_1\times G_2$ such that
$\mathbf{E}(I_H(\gamma_1,\gamma_2))=1$ and for every $(g_1,g_2)\in H$
\[
\Xi_{12}^K(g_1,g_2,A) = \int I_A(k) h(g_1,g_2,k)
\Xi_{1}^K(g_1,dk)\qquad
\mbox{for all }A\in\mathcal{K}.
\]
Then there is $H'\subset G_1\times K$ with
$\mathbf{E}(I_{H'}(\gamma_1,\kappa))=1$ so that for all $(g_1,k)\in H'$
\[
\Xi_{1K}^2(g_1,k,B) = \int I_B(g_2) h(g_1,g_2,k)
\Xi_{1}^2(g_1,dg_2)\qquad
\mbox{for all }B\in\mathcal{G}_2.
\]
\end{lem}

We now recall two results of von Weizs{\"a}cker that are of central
importance in our proofs. The first result is a special case of the
result in~\cite{vW83}, pages 95 and 96.
%
%le8.3 #&#
\begin{lem}
\label{lemvWconditioning}
Let $G$, $G'$ and $H$ be Polish spaces, and denote by $g$, $g'$ and~$h$
the canonical projections from $G\times G'\times H$ on $G$, $G'$ and $H$,
respectively. Let $\mathbf{Q}$ be a probability measure on $G\times
G'\times H$, and let $q_{\cdot,\cdot}\dvtx G\times
G'\times\mathcal{B}(H)\to[0,1]$ and $q_{\cdot}\dvtx G\times\mathcal
{B}(G'\times
H)\to[0,1]$ be versions of the regular conditional probabilities
$\mathbf{Q}[h\in\cdv g,g']$ and
$\mathbf{Q}[(g',h)\in\cdv g]$, respectively.
Then for $\mathbf{Q}$-a.e. $x\in G$, the kernel
$q_{x,g'}[ \cdot]$ is a version of the regular conditional probability
$q_{x}[h\in\cdv g']$.
\end{lem}

Though the second result is not given precisely in this form in
\cite{vW83}, its proof follows easily from~\cite{vW83} modulo minor
modifications (see also~\cite{vH10}, Section 4.1).
%
%le8.4 #&#
\begin{lem}
\label{lemvWexchg}
Let $G$ and $H$ be Polish spaces, let $(X_n)_{n\ge0}$ be a sequence
of random variables with values in $G$ and let $Y$ be a random
variable with values in $H$ on some underlying probability space
$(\Omega,\mathcal{F},\mathbf{P})$. Define the $\sigma$-field
$\mathcal{H}=\sigma\{Y\}$ and the decreasing filtration
$\mathcal{G}_n=\sigma\{X_k\dvtx k\ge n\}$. Then
\[
\bigcap_{n\ge0} \mathcal{H}\vee\mathcal{G}_n =
\mathcal{H},\qquad
\mathbf{P}\mbox{-a.s.}
\]
if and only if
\[
\bigcap_{n\ge0}\mathcal{G}_n
\mbox{ is }
\mathbf{P}^{\mathcal{H}}\mbox{-trivial},\qquad
\mathbf{P}\mbox{-a.s.},
\]
where $\mathbf{P}^{\mathcal{H}}$ is a version of the
regular conditional probability
$\mathbf{P}((X_n)_{n\ge0}\in\cdv \mathcal{H})$.\vadjust{\goodbreak}
\end{lem}

%s9 #&#
\section{Notation list}
\label{secnotationlist}

The following list of frequently used notation, together with the page
numbers where they are defined, is included for easy reference.\vspace*{6pt}%\newpage

\begin{center}
\begin{tabular*}{\tablewidth}{@{}lp{323pt}@{}}
 $E$\hphantom{$(G)$} &  State space of unobservable component $X_n$\dotfill
                \nompageref{page7}\\[3pt]
 $F$ &  State space of observable component $Y_n$\dotfill
                \nompageref{page7}\\[3pt]
 $P$ &  Transition kernel of $(X_n,Y_n)_{n\in\mathbb{Z}}$\dotfill
                \nompageref{page7a}\\[3pt]
 $P'$ &  Transition kernel of the reversed model $(X_{-n},Y_{-n})_{n\in\mathbb{Z}}$\dotfill
                \nompageref{page7aa}\\[3pt]
 $P^X$ &  Conditional transition kernel of $(X_n)_{n\in\mathbb{Z}}$ given $(Y_n)_{n\in\mathbb{Z}}$\dotfill
                \nompageref{page24}\\[3pt]
 $P_0$ &  Reference kernel on $E$ such that $P\sim P_0\otimes Q$ (Assumption~\ref{asptnondeg})\dotfill
                \nompageref{page9}\\[3pt]
 $P_n^X$ &  Version of $\mathbf{P}(X_n\in \cdv X_0)$ (Lemma~\ref{lemnondegexact})\dotfill
                \nompageref{page16}\\[3pt]
 $P_n^Y$ &  Version of $\mathbf{P}(Y_n\in \cdv Y_0)$ (Lemma~\ref{lemnondegexact})\dotfill
                \nompageref{page16}\\[3pt]
 $Q$ &  Reference kernel on $F$ such that $P\sim P_0\otimes Q$ (Assumption~\ref{asptnondeg})\dotfill
                \nompageref{page9}\\[3pt]
 $U$ &  Filter recursion (Lemma~\ref{lemintrofilt}; cf. Section~\ref{secpairm})\dotfill
                \nompageref{page8}\\[3pt]
 $X_n$ &  Unobservable component of model\dotfill
                \nompageref{page7aaa}\\[3pt]
 $Y$ &  Observation path $(Y_k)_{k\in\mathbb{Z}}$\dotfill
                \nompageref{page7b}\\[3pt]
 $Y_n$ &  Observable component of model\dotfill
                \nompageref{page7aaa}\\[3pt]
 $\Gamma$ &  Transition kernel of $(\Pi_n^\mu,Y_n)_{n\ge 0}$ (Lemma~\ref{lemintromarkov}; cf. Section~\ref{secpairm})\dotfill
                \nompageref{page8a}\\[3pt]
 $\Lambda$ &  Transition kernel of $(\Pi_n^\mu,X_n,Y_n)_{n\ge 0}$ (Lemma~\ref{lemintromarkov}; cf. Section~\ref{sectriplem})\dotfill
                \nompageref{page8a}\\[3pt]
 $\Omega$ &  Canonical path space\dotfill
                \nompageref{page7bb}\\[3pt]
 $\Omega^X$ &  Canonical path space of unobservable component\dotfill
                \nompageref{page7bb}\\[3pt]
 $\Omega^Y$ &  Canonical path space of observable component\dotfill
                \nompageref{page7bb}\\[3pt]
 $\Pi_n^\mu$ &  The nonlinear filter $\mathbf{P}^\mu(X_n\in \cdv \mathcal{F}^Y_{[0,n]})$ (cf. Lemma~\ref{lemrefprob})\dotfill
                \nompageref{pagea8}\\[3pt]
 $\Theta$ &  The canonical shift on $\Omega$\dotfill
                \nompageref{page7b}\\[3pt]
 $\mathbf{P}$ &  Law of $(X_n,Y_n)_{n\in\mathbb{Z}}$\dotfill
                \nompageref{page7a}\\[3pt]
 $\mathbf{P}^Y$ &  Law of the observations $(Y_n)_{n\in\mathbb{Z}}$\dotfill
                \nompageref{page24}\\[3pt]
 $\mathbf{P}^\mu$ &  $\int \mathbf{P}^{z,w} \mu(dz,dw)$\dotfill
                \nompageref{pagePP}\\[3pt]
 $\mathbf{P}^{z,w}$ &  Conditional law of $(X_n,Y_n)_{n\in\mathbb{Z}}$ given $X_0=z$, $Y_0=w$\dotfill
                \nompageref{pageCond}\\[3pt]
 $\mathbf{P}_{z,y}$ &  Conditional law of $(X_n)_{n\ge 0}$ given $X_0=z$, $Y=y$\dotfill
                \nompageref{pageCondlaw}\\[3pt]
 $\mathcal{B}(G)$ &  Borel $\sigma$-field of $G$\dotfill
                \nompageref{pageBorel}\\[3pt]
 $\mathcal{F}$ &  Borel $\sigma$-field of $\Omega$\dotfill
                \nompageref{pageBorelsigma}\\[3pt]
 $\mathcal{F}^Z$ &  $\mathcal{F}_{\mathbb{Z}}^Z$\dotfill
                \nompageref{page7b}\\[3pt]
 $\mathcal{F}^Z_+$ &  $\mathcal{F}_{[0,\infty[}^Z$\dotfill
                \nompageref{page7b}\\[3pt]
 $\mathcal{F}_I$ &  $\mathcal{F}_I^X\vee\mathcal{F}_I^Y$\dotfill
                \nompageref{pageFF}\\[3pt]
 $\mathcal{F}_I^Z$ &  $\sigma\{Z_k\dvtx k\in I\}$\dotfill
                \nompageref{pageFF}\\[3pt]
\end{tabular*}
\end{center}

\begin{center}
\begin{tabular*}{\tablewidth}{@{}lp{323pt}@{}}
 $\mathcal{F}_n$ &  $\mathcal{F}_n^X\vee\mathcal{F}_n^Y$\dotfill
                \nompageref{pageFFY}\\[3pt]
 $\mathcal{F}_n^Z$ &  $\mathcal{F}_{]-\infty,n]}^Z$\dotfill
                \nompageref{pageFFY}\\[3pt]
 $\mathcal{P}(G)$ &  Space of probability measures on $G$\dotfill
                \nompageref{pageBorel}\\[3pt]
 $\mu_w$ &  Version of $\mathbf{P}^\mu(X_0\in \cdv Y_0)$\dotfill
                \nompageref{page30}\\[3pt]
 $\pi$ &  Invariant measure of $(X_n,Y_n)_{n\in\mathbb{Z}}$\dotfill
                \nompageref{page7a}\\[3pt]
 $\pi^X$ &  Version of $\mathbf{P}(Y_0\in \cdv X_0)$ (Lemma~\ref{lemnondegexact})\dotfill
                \nompageref{page16}\\[3pt]
 $\pi^Y$ &  Version of $\mathbf{P}(X_0\in \cdv Y_0)$ (Lemma~\ref{lemnondegexact})\dotfill
                \nompageref{page16}\\[3pt]
 $\varpi$ &  Conditional law of $X_0$ given $(Y_n)_{n\in\mathbb{Z}}$\dotfill
                \nompageref{page24}\\[3pt]
 $\varpi^+$ &  Conditional law of $X_0$ given $(Y_n)_{n\ge 0}$\dotfill
                \nompageref{page28}\\[3pt]
 $b\mathsf{m}$ &  Barycenter of $\mathsf{m}$\dotfill
                \nompageref{pageBary}\\[3pt]
 $g$ &  Transition density of $P$ (Assumption~\ref{asptnondeg})\dotfill
                \nompageref{page9}\\[3pt]
\end{tabular*}
\end{center}
\end{appendix}

%suskaldyti doi

% imsref loaded by lrinkeviciute, 2012-02-29 09:06:49

\printaddresses

\end{document}